\documentclass[reqno]{amsart}
\usepackage[normalem]{ulem}

\usepackage[utf8x]{inputenc}
\usepackage[T1]{fontenc}

\usepackage{enumerate}

\usepackage{subfig}
\usepackage[pdftex]{graphicx} 
\usepackage[pdftex,linkcolor=black,pdfborder={0 0 0}]{hyperref} 
\usepackage{calc} 

\usepackage[a4paper, lmargin=0.1666\paperwidth, rmargin=0.1666\paperwidth, tmargin=0.1111\paperheight, bmargin=0.1111\paperheight]{geometry} 

\graphicspath{{Figures/}}

\usepackage[all]{nowidow} 
\usepackage[protrusion=true,expansion=true]{microtype} 

\usepackage{amsmath, amssymb, amsthm}
\usepackage{setspace}
\usepackage{mathrsfs}
\usepackage{graphicx}
\usepackage{pdfpages}
\usepackage{listings}
\usepackage{bbm}
\usepackage{hyperref}
\usepackage{subfig}
\usepackage{color}
\usepackage{float}
\usepackage{amssymb}

\usepackage{cancel}

\DeclareMathOperator{\BPP}{\mathbb{P}}

\DeclareMathOperator{\CTT}{\mathcal{T}}

\DeclareMathOperator{\La}{\Lambda}

\newcommand{\rr}{\mathbb{R}}
\newcommand{\pp}{\mathbb{P}}
\newcommand{\E}{\mathbb{E}}
\newcommand{\teps}{\epsilon}

\theoremstyle{plain}
\newtheorem{thm}{Theorem}[section]
\newtheorem{dfn}[thm]{Definition}
\newtheorem{ass}[thm]{Assumption}
\newtheorem{prop}[thm]{Proposition}
\newtheorem*{prop*}{Proposition}
\newtheorem{lem}[thm]{Lemma}
\newtheorem{cor}[thm]{Corollary}

\theoremstyle{definition}
\newtheorem{rmk}[thm]{Remark}

\begin{document}

\title[Time-Stepping for a McKean-Vlasov Equation with Blow-Ups]{Convergence of a time-stepping scheme to the free boundary in the supercooled Stefan problem}

\author[V.~Kaushansky]{Vadim Kaushansky}
\address{Department of Mathematics, University of California, Los Angeles, CA 90095, U.S.}
\email{vadim.kaushansky@gmail.com}
\author[C.~Reisinger]{Christoph Reisinger}
\address{Mathematical Institute, University of Oxford, Andrew Wiles Building, Radcliffe Observatory Quarter, OX2 6GG, Oxford, U.K.}
\email{christoph.reisinger@maths.ox.ac.uk}
\author[M.~Shkolnikov]{Mykhaylo Shkolnikov}
\address{ORFE Department, Bendheim Center for Finance, and Program in Applied \& Computational Mathematics, Princeton University, Princeton, NJ 08544, U.S.}
\email{mshkolni@gmail.com}
\author[Z.Q.~Song]{Zhuo Qun Song}
\address{Department of Mathematics, Princeton University, Princeton, NJ 08544, U.S.}
\email{zhuoqunsong@gmail.com}



\begin{abstract}
The supercooled Stefan problem and its variants describe the freezing of a supercooled liquid in physics, as well as the large system limits of systemic risk models in finance and of integrate-and-fire models in neuroscience. Adopting the physics terminology, the supercooled Stefan problem is known to feature a finite-time blow-up of the freezing rate for a wide range of initial temperature distributions in the liquid. Such a blow-up can result in a discontinuity of the liquid-solid boundary. In this paper, we prove that the natural Euler time-stepping scheme applied to a probabilistic formulation of the supercooled Stefan problem converges to the liquid-solid boundary of its physical solution globally in time, in the Skorokhod M1 topology. In the course of the proof, we give an explicit bound on the rate of local convergence for the time-stepping scheme. We also run numerical tests to compare our theoretical results to the practically observed convergence behavior.
\end{abstract}

\maketitle


\section{Introduction}

The classical formulation of the one-dimensional one-phase Stefan problem (henceforth, simply Stefan problem), introduced by \textsc{Stefan} in \cite{Stefan1,Stefan2,Stefan3,Stefan4} (see also \cite{LC}), can be stated as follows: 
\begin{equation}\label{eq:Stefan}
	\begin{split}
		& \partial_t u = \frac{1}{2}\partial_{xx}u,\;\; x\ge \Lambda_t,\;\; t\ge 0, \\
		& u(0,x)=f(x),\;\; x\ge0 \quad\text{and}\quad u(t,\Lambda_t)=0,\;\; t\ge0, \\
		& \dot{\Lambda}_t = \frac{\alpha}{2}\partial_x u(t,\Lambda_t),\;\; t\ge0\quad\text{and}\quad\Lambda_0=0.
	\end{split}
\end{equation}
Hereby, the negative $-f$ of the given function $f:\,[0,\infty)\to\rr$ stands for the initial temperature distribution in a liquid relative to its equilibrium freezing point; $\Lambda_t$ is the unknown location of the liquid-solid boundary at time $t$; the negative $-u(t,\cdot)$ of the unknown function $u(t,\cdot):\,[\Lambda_t,\infty)\to\rr$ represents the temperature distribution in the liquid relative to its equilibrium freezing point at time $t$; and $\frac{2}{\alpha}>0$ is the density of latent heat when the heat capacity of the liquid is normalized to $2$ and its thermal conductivity to $1$. We consider the Stefan problem \eqref{eq:Stefan} in the \textit{supercooled} regime, i.e., when $f\ge0$.  

\medskip

As first observed by \textsc{Sherman} in \cite{Sher}, classical solutions of the supercooled Stefan problem can exhibit a finite-time blow-up, in the sense that $\lim_{t\uparrow t_*} \dot{\Lambda}_t=\infty$ for some $t_*\in(0,\infty)$. Subsequently, much attention has been devoted to the construction of classical solutions to \eqref{eq:Stefan} on $[0,t_*)$, to the distinction between the cases $t_*<\infty$ and $t_*=\infty$ in terms of the initial condition $f$ (see \cite{FP1,FP2,DF,HOL,LO,FPHO} and the references therein), and to the regularization of \eqref{eq:Stefan} when $t_*<\infty$ (see \cite{Vis,DHOX,FPHO,HX,Wei} and the references therein). 

\medskip

More recently, it has been discovered in \cite{NS1} that (a variant of) \eqref{eq:Stefan} admits a probabilistic reformulation (see also \cite[Section 1.2]{DIRT2} for a related discussion). Indeed, suppose that $\|f\|_{L^1([0,\infty))}=1$, let $X_{0-}$ be a random variable with the probability density $f$, and consider the problem of finding a non-decreasing function $\Lambda$ such that the stochastic process 
\begin{equation}\label{eq::init1}
	X_t = X_{0-} + B_t - \Lambda_t,\;\; t\ge 0 
\end{equation}
satisfies the constraint
\begin{equation}\label{eq::init2}
	\Lambda_t = \alpha\pp\Big(\inf_{0\le s\le t} X_s \le 0\Big),\;\;t\ge 0, 
\end{equation}
where $B$ is a standard Brownian motion independent of $X_{0-}$. Assuming that $f$ belongs to the Sobolev space $W^1_2([0,\infty))$ and that $f(0)=0$,  there exists a $T\in(0,\infty)$ and a solution $\Lambda$ of \eqref{eq::init1}--\eqref{eq::init2} on $[0,T]$ with $\dot{\Lambda}\in L^2([0,T])$; moreover,  with $\tau:=\inf\{t\ge0:\,X_t\le 0\}$, the densities $p(t,\cdot)$, $t\in[0,T]$ of $X_t\,\mathbf{1}_{\{\tau>t\}}$, $t\in[0,T]$ on $(0,\infty)$ form the unique solution of
\begin{equation}\label{p_PDE}
	\begin{split}
		& \partial_t p=\frac{1}{2}\partial_{xx}p+\dot{\Lambda}_t\partial_x p,\;\;x\ge0,\;\;t\in[0,T], \\ 
	 	& p(0,x)=f(x),\;\; x\ge 0\quad\text{and}\quad p(t,0)=0,\;\;t\in[0,T], \\
		& \dot{\Lambda}_t=\frac{\alpha}{2}\partial_x p(t,0),\;\;t\in[0,T]\quad\text{and}\quad\Lambda_0=0
	\end{split}
\end{equation}
in the Sobolev space $W^{1,2}_2([0,T]\times[0,\infty))$ (cf.~\cite[Propositions 4.1, 4.2]{NS1}), i.e., $u(t,x):=p(t,x-\Lambda_t)$, $x\ge\Lambda_t$, $t\in[0,T]$ is the solution of \eqref{eq:Stefan} in $W^{1,2}_2(\{(t,x)\in[0,T]\times[0,\infty):\,x\ge\Lambda_t\})$.

\medskip

Due to the absence of the derivative $\dot{\Lambda}$ in the probabilistic formulation \eqref{eq::init1}--\eqref{eq::init2}, the latter allows to study the liquid-solid boundary $\Lambda$ for all times. As noted in \cite[Theorem 1.1]{HLS}, $\Lambda$ cannot be continuous for all $t\in[0,\infty)$ when $\E[X_{0-}]<\frac{\alpha}{2}$, which necessitates to consider the solutions of \eqref{eq::init1}--\eqref{eq::init2} in the space $D([0,\infty))$ of right-continuous functions with left limits in general. However, in $D([0,\infty))$ uniqueness does not hold for \eqref{eq::init1}--\eqref{eq::init2}, since the jump sizes of $\Lambda$ are not determined uniquely by \eqref{eq::init1} and \eqref{eq::init2} alone (see, e.g., \cite[Example 2.2, Figure 3]{HLS}, \\ as well as \cite[p.~7, last paragraph]{NS2}, \cite[discussion preceding Definition 2.2]{DIRT}). 
We will consider those solutions where the jump sizes are chosen as
\begin{equation}\label{physical}
	\Lambda_t-\Lambda_{t-}=\inf\Big\{x>0:\,\pp\big(\tau\ge t,\,X_{t-}\in(0,x]\big)<\frac{x}{\alpha}\Big\},\;\;t\ge 0, 
\end{equation} 
with $\Lambda_{t-}:=\lim_{s\uparrow t} \Lambda_s$ and $X_{t-}:=\lim_{s\uparrow t} X_s$. Physically, \eqref{physical} states that the supercooled liquid heats up to its equilibrium freezing temperature and freezes on the smallest interval $[\Lambda_{t-},\Lambda_{t-}+x)$ for which this transition is energy neutral. Similarly, when \eqref{eq::init1}--\eqref{eq::init2} and its variants arise in the large system limit of systemic risk models in finance (see \cite{NS1,HLS,LS}), or of integrate-and-fire models in neuroscience (see \cite{DIRT}), only the solutions with the minimality property \eqref{physical} seem economically, or biologically, relevant. Under the Assumption \ref{ass}(a), the solution $\Lambda$ of \eqref{eq::init1}--\eqref{eq::init2} in $D([0,\infty))$ satisfying \eqref{physical} is unique (see \cite[Theorem 1.4]{DNS}), and it is referred to as the \textit{physical solution} of \eqref{eq::init1}--\eqref{eq::init2}. 

\begin{ass}\label{ass}
\begin{enumerate}[(a)]
	\item $X_{0-}\!\!\ge\!0$ possesses a bounded density $f$ on $[0, \infty)$ that changes monotonicity finitely often on compacts.
	\item There exists a physical solution of \eqref{eq::init1}--\eqref{eq::init2} starting from $X_{0-}$. 
\end{enumerate}	
\end{ass}

\noindent The existence of a physical solution is shown for $X_{0-}\!\ge\!0$ with $\E[X_{0-}]\!<\!\infty$ in \cite[Theorem~6.5]{CRS}, building in part on \cite[proof of Theorem 3.2]{LS} which made stronger assumptions on $X_{0-}$ (see also \cite[Theorem 4.4]{DIRT2}, \cite[Theorem 2.3]{NS1} for existence results in related contexts).

\medskip

The absence of a closed-form solution calls for numerical methods. Proposed methods include approximations to Volterra equations that $\Lambda$ solves (see \cite{lipton2019semi}) and a fixed-point iteration in the space of boundaries (see \cite{CRS}). Herein, we follow the general strategy of \cite{KR}:
\begin{itemize}
	\item Define a discrete process $\widetilde{X}$ and a discrete function $\widetilde{\Lambda}$ such that $\widetilde{X}$ and $\widetilde{\Lambda}$ are close to $X$ and $\Lambda$ in some sense.
	\item Approximate $\widetilde{X}$ and $\widetilde{\Lambda}$ by a representative particle in a particle system $\{\overline{X}^i\}_{i=1}^N$ and a corresponding estimator $\overline{\Lambda}^N$, which can be generated through simulation.
\end{itemize}
Subsequently, given a discrete time step $\Delta>0$, the following natural definition of $(\widetilde{X},\widetilde{\Lambda})$, referred to as $(X^\Delta,\Lambda^\Delta)$, is chosen in \cite{KR}.   

\begin{dfn}\label{defn::discrete}
	Let $X^\Delta_{n\Delta}$, $n=0,\,1,\,\ldots$ and $\Lambda^\Delta_{n\Delta}$, $n=1,\,2,\,\ldots$ be defined recursively by $X^\Delta_0=X_0$, 
	\begin{equation}\label{eq::init_discrete}
		\Lambda^{\Delta}_{n\Delta}=\alpha\pp\Big(\min_{0\le m\le n-1} X^\Delta_{m\Delta}<0\Big),\;n=1,2,\ldots\;\text{and}\; X^\Delta_{n\Delta}=X_{0-}+B_{n\Delta}-\Lambda^\Delta_{n\Delta},\;n=1,2,\ldots.  
	\end{equation}
	We extrapolate by setting $\Lambda^\Delta_{0-}=0$, $\Lambda^\Delta_0=\alpha\pp(X_0\le0)$ and $\Lambda^{\Delta}_t=\Lambda^{\Delta}_{\lfloor t/\Delta \rfloor \Delta}$, $t>0$, as well as $X^\Delta_{0-}=X_{0-}$ and $X^{\Delta}_t = X_{0-} + B_t - \Lambda^{\Delta}_t$, $t\ge0$. 
\end{dfn}

It is shown in \cite[Theorem 1]{KR} that under the additional decay assumption $f(x)=O(x^\beta)$, $x\downarrow0$ for some $\beta\in(0,1]$ (originally introduced in \cite{HLS}) the convergence $\Lambda^\Delta\underset{\Delta\downarrow0}{\longrightarrow}\Lambda$ holds uniformly in time on a sufficiently small interval $[0,T]$, depending on the model parameters and preceding the first jump of $\Lambda$, with an order arbitrarily close to $\frac{1}{2}$. In the classical setting of stochastic differential equations (SDEs) with Lipschitz coefficients, the strong order of the Euler-Maruyama scheme, of which the above is a variant, is $\frac{1}{2}$, but increases to $1$ if the diffusion coefficient is constant (as is the case here) since the Euler scheme coincides with the first order Milstein scheme in this case. 

The reduction of the convergence rate in the present setting comes from (a) the possible non-Lipschitzianity of the densities of $X_t\,\mathbf{1}_{\{\tau>t\}}$ at $0$ and (hence) of the bounded variation part $-\Lambda$ of $X$, and (b) the underestimation of the hitting probability between time steps. In the setting of \cite{KR}, these difficulties can be resolved to obtain convergence of order $1$ in the following way (see \cite[Corollary 1]{KR}): refining the time mesh close to time $0$ where $\dot{\Lambda}$ has a singularity, and
using a continuous-time Brownian bridge interpolation of the Euler--Maruyama scheme as in \cite{gobet2000weak} to correct the leading order term in $\Delta$ of the hitting probability.
As the loss function with interpolation lies between the approximate loss function for our time-stepping scheme and the true loss function, the global convergence result below extends automatically to the scheme with interpolation. An improvement of the convergence order will only be observed locally where the loss function is sufficiently regular. The time mesh refinement in  \cite{KR} relies on explicit knowledge of the polynomial order of the singularity of $\dot{\Lambda}$ at $t=0$.
Although we give a modulus of right-continuity for $\Lambda$ in Corollary~\ref{cor:mod_of_cont}, it is unclear how it can be leveraged to control the local error by mesh refinement and obtain a higher convergence order. It may be possible to construct an adaptive scheme which estimates heuristically the local regularity of $\Lambda$, in particular in the run-up to a singularity, but this is beyond the scope of the present paper.

\medskip

Our main result is concerned with the convergence of $\Lambda^\Delta$ to $\Lambda$ \textit{globally} in time. For each $T\in(0,\infty)$, we let $D([0,T])$ be the space of real-valued functions on $[0,T]$ that are right-continuous at all $t\in[0,T)$ and have left limits at all $t\in(0,T]$. We endow $D([0,T])$ with the topology of Skorokhod M1 convergence from \cite{Sko}, whose definition for non-decreasing functions is recalled in Subsection \ref{subsec:M1} for the convenience of the reader. Our main result can then be stated as follows. 

\begin{thm}\label{thm::main}
Under Assumption \ref{ass} and for any continuity point $T\in(0,\infty)$ of $\Lambda$, it holds $\Lambda^\Delta|_{[0,T]}\overset{\emph{M1}}{\underset{\Delta\downarrow0}{\longrightarrow}}\Lambda|_{[0,T]}$ in $D([0,T])$.   
\end{thm}

\begin{rmk}
	The notion of convergence in Theorem \ref{thm::main} is often referred to as M1 convergence in $D([0,\infty))$, see, e.g., \cite[Section 12.9]{Whi}. In particular, for any $T\in(0,\infty)$ such that $\Lambda|_{[0,T]}$ is continuous, this notion of convergence yields the uniform convergence $\sup_{t\in[0,T]} |\Lambda^\Delta_t-\Lambda_t|\underset{\Delta\downarrow0}{\longrightarrow}0$ (use, e.g., \cite[Corollary 12.5.1, then implication (v)$\Rightarrow$(i) in Theorem 12.4.1, and finally Lemma 12.4.2]{Whi}). Thus, Theorem \ref{thm::main} naturally includes the convergences addressed in \cite[Theorems 1, 4]{KR}. Moreover, Theorem \ref{thm::main} gives the almost sure convergence of $X^\Delta=X_{0-}+B-\Lambda^\Delta$ to $X=X_{0-}+B-\Lambda$ with respect to the M1 convergence in $D([0,\infty))$, thanks to \cite[Corollary 12.7.1]{Whi}.
	
\end{rmk}

A crucial and novel element of our analysis is the argument that for vanishing mesh size the numerical solution approaches the true physical solution immediately after a jump. This is achieved by first showing that the probability density of the numerical solution is arbitrarily close to that of the true physical solution immediately prior to the jump, and then using this property to prove that the numerical solution ``catches up'' arbitrarily fast with the physical solution after the jump of the latter (see Figure \ref{fig:example2} for an illustration). The convergence of the explicit Euler-type scheme to the \emph{physical} solution herein is akin to that of a system with delays to the physical solution in \cite[Theorem 4.9]{DIRT}. The physicality of the limit can be rationalized by noticing that the numerical approximations converge from below, and reveals that other than by a minimality property the physical solution is also characterized as the continuous-time limit of a discrete approximation with a delayed self-reinforcement effect. This differs from the proofs of \cite[Theorem~4.7]{DIRT}, \cite[Theorem 2.3]{NS1}, \cite[Theorem 3.2]{LS} which show the convergence of approximating particle systems to the physical solution. There, it is hard to rule out ``overshooting'' (non-physicality) of the limit.

\medskip

In the course of the proof of Theorem \ref{thm::main} we obtain an explicit bound on the rate of convergence \textit{locally} in time, which varies according to the density of $X_0$ on a right neighborhood of $0$. This result, interesting on its own, is reported in Theorem \ref{thm:loc_conv}. 

\begin{thm}\label{thm:loc_conv}
Suppose Assumption \ref{ass}. Let the density of $X_0$ be less or equal to $\frac{1}{\alpha}-\psi$ on an interval $(0,\delta]$, with a strictly increasing function $\psi:\,(0,\delta]\to(0,\infty)$. For $x\in(0,\delta]$, define $\Psi(x)\!=\!\int_0^x \psi(y)\,\mathrm{d}y$ and
\begin{equation}
\widetilde{\Psi}(x)=\Bigg\lbrace\!\!\begin{array}{ll} 
c_1x & \text{if}\;\;\psi(0):=\lim_{y\downarrow 0} \psi(y)>0, \\
c_2\int_0^x \Phi\big(q-\frac{c_3}{\Psi(y)}\big)\,\mathrm{d}y & \text{if}\;\;\psi(0)=0,
\end{array}
\end{equation}
where $\Phi$ is the standard Gaussian cumulative distribution function, 
 \begin{eqnarray}
 && c_1:=\frac{1}{4}\Big(\psi(0)-\big(2\|f\|_\infty+\psi(0)-1/\alpha\big)_+\,e^{-\frac{\delta^2}{3\teps}}\Big)>0, \qquad \\
&& c_2:=\psi(\delta/2)-\big(2\|f\|_\infty+\psi(\delta/2)-1/\alpha\big)_+\,e^{-\frac{7\delta^2}{24\teps}}>0, \qquad \\
&& c_3:=\big(\|f\|_\infty\,\sqrt{2/\pi}-q\psi(\delta)\big)\,\delta/2>0, \qquad \\
&& q:=\Phi^{-1}(1/4)<0, \qquad \\
&& 0<\teps<\begin{cases}
\min\Big(\frac{\pi}{8\|f\|_\infty^2}\,\Psi(\delta/6)^2,\,\frac{\delta^2}{3\log\max(1,(2\|f\|_\infty+\psi(0)-1/\alpha)_+/\psi(0))}\Big) & \text{if}\;\;\psi(0)\!>\!0, \\
\frac{\pi}{8\|f\|_\infty^2}\,\Psi(\delta/6)^2 & \text{if}\;\;\psi(0)\!=\!0.
\end{cases} \label{teps:UBD} \qquad
\end{eqnarray} 	
Then, for all $\teps$ adhering to \eqref{teps:UBD} and all $\Delta>0$ small enough,
\begin{eqnarray}
\sup_{s\in[0,\teps]} |\Lambda^\Delta_s\!-\!\Lambda_s| \le
\widetilde{\Psi}^{-1}\Big(\|f\|_\infty\Big(\!192\sqrt{2\Delta\log(\lceil\teps/\Delta\rceil/2)} 
\!+\!2\Psi^{-1}\big(2\sqrt{\pi}^{-1}\|f\|_\infty\sqrt{\Delta}\big)\!\Big)\!\Big) \qquad\; && \label{roc} \\
+\max(\alpha\|f\|_\infty,1)\Big(\!192\sqrt{2\Delta\log(\lceil\teps/\Delta\rceil/2)} 
\!+\!2\Psi^{-1}\big(2\sqrt{\pi}^{-1}\|f\|_\infty\sqrt{\Delta}\big)\!\Big). && \nonumber
\end{eqnarray}
\end{thm}

\begin{rmk}
Assumption \ref{ass}(a) and the condition \eqref{physical} together imply that the density of $X_0$ is necessarily less or equal to $\frac{1}{\alpha}-\psi$ on an interval $(0,\delta]$, for some strictly increasing function $\psi:\,(0,\delta]\to(0,\infty)$. If the density of $X_0$ is such that $\delta>0$ can be chosen arbitrarily large, then the upper bound on $\teps$ in \eqref{teps:UBD} can be made arbitrarily large, rendering \eqref{roc} a \textit{global} estimate. In general, however, to control the rate of convergence globally one needs a quantitative bound on the density of $X_t\,\mathbf{1}_{\{\tau>t\}}$ on a right neighborhood of $0$ for an arbitrary $t\ge0$, which has proved elusive. 
\end{rmk}

\begin{rmk}
\label{rm_order}
		\!\!(a) When $\psi(0)\!>0$, we may swap $\psi$ for the constant $\psi(0)$, simplifying \eqref{roc} to 
			\begin{equation}\label{roc_simplified}
			\sup_{s\in[0,\teps]} |\Lambda^\Delta_s-\Lambda_s| \le
			c_4\sqrt{\Delta}\,\big(\sqrt{\log(\lceil\teps/\Delta\rceil/2)}+c_5\big)
			\end{equation}
		and \eqref{teps:UBD} to $0<\teps<c_6\delta^2$, with 
		$c_4:=192\sqrt{2}\,\big(\|f\|_\infty/c_1+\max(\alpha\|f\|_\infty,1)\big)$, 
		$c_5:=\frac{\|f\|_\infty}{48\sqrt{2\pi}\psi(0)}$, and 
		$c_6:=\min\big(\frac{\pi\psi(0)^2}{288\|f\|_\infty^2},\,\frac{1}{3\log\max(1,(2\|f\|_\infty+\psi(0)-1/\alpha)_+/\psi(0))}\big)$. 
		As we detail in Subsection \ref{sim1/2}, the theoretical order $\frac{1}{2}$ in \eqref{roc_simplified} agrees with the one found in numerical simulations. 
		
		\vspace{5pt}
		
		\noindent (b) When $\psi(0)=0$, we have $\lim_{x\downarrow0} \widetilde{\Psi}^{-1}(x)=0$ and $\lim_{x\downarrow0} (\widetilde{\Psi}^{-1})'(x)=\infty$, so that the $\widetilde{\Psi}^{-1}$-term on the right-hand side of \eqref{roc} determines its order of magnitude in $\Delta$. As we illustrate in 
		Subsection \ref{sec5.3}, the order observed in numerical simulations agrees instead with that of the term on the second line in \eqref{roc}. This leads us to conjecture that, in fact, the latter gives the true order of the left-hand side of \eqref{roc} in $\Delta$, whereas the $\widetilde{\Psi}^{-1}$-term on the right-hand side of \eqref{roc} is an artifact of our proof technique. Specifically, if $\partial_x^b \psi(0+) = 0$, $b=0,\,1,\,\ldots,\,a-1$, $\partial_x^a \psi(0+) > 0$, as for the monomials $x^a$ in Subsection \ref{sec5.3},
it would yield the order $1/(2(a+1))$ noticed in the numerical experiments. For $\psi$ with vanishing derivatives at $0$ of all orders, we do not expect to see a positive covergence order, but it is hard to check this claim numerically.
\end{rmk}

The rest of the paper is structured as follows. In Section \ref{sec:prelim}, we collect a number of preliminaries: a classification of time points for the physical solution (Subsection \ref{classi}); explicit estimates on the densities of $X_s\,\mathbf{1}_{\{\tau>s\}}$ and $\min_{0\le m\le n} X_{t+m\Delta}\,\mathbf{1}_{\{\tau>t\}}$ near $0$, as well as on the increments of $\Lambda$, in right time neighborhoods (Subsection \ref{rightn}); and definitions of the Skorokhod M1 convergence for non-decreasing functions and of some auxiliary time-stepping schemes (Subsection \ref{subsec:M1}). In Section \ref{sec:loc_conv}, we prove Theorem \ref{thm:loc_conv} by an inductive argument relying on the results of Subsection \ref{rightn}. Subsequently, in Section \ref{sec:mainproof}, we show our main result (Theorem \ref{thm::main}) by extending Theorem \ref{thm:loc_conv} to the statement that the M1 convergence $\Lambda^\Delta\underset{\Delta\downarrow0}{\longrightarrow}\Lambda$ on an interval $[0,T)$ implies the M1 convergence $\Lambda^\Delta|_{[0,T+\teps]}\underset{\Delta\downarrow0}{\longrightarrow}\Lambda|_{[0,T+\teps]}$ for some $\teps>0$. We achieve this separately for the cases that $T$ is a continuity point of $\Lambda$ (Subsection~\ref{sec::no_jump}) and that $T$ is a discontinuity point of $\Lambda$ (Subsection \ref{sec::special_case}). In both cases, we build on Theorem \ref{thm:loc_conv} and the error propagation bounds of Subsection \ref{sec::stab}. Subsection \ref{sec::post_jump} then offers the main line of the proof of Theorem \ref{thm::main}. In Section \ref{sec::numerical}, we give the results of our simulations, based on a particle approximation, and compare them with the theoretical findings. The convergence rate of this particle approximation is controlled in the appendix.

\medskip

\noindent\textbf{Acknowledgement.} M.~Shkolnikov has been partially supported by the NSF grants DMS-1811723, DMS-2108680 and a Princeton SEAS innovation research grant.


\section{Preliminaries} \label{sec:prelim}

\subsection{Classification of time points}\label{classi}

The basis of our analysis is the following proposition (see \cite[Theorem 1.4, first statement of Theorem 1.1]{DNS} and recall the minimality property~\eqref{physical}).
\begin{prop}\label{thm1}
Suppose Assumption \ref{ass}. Then, the physical solution is unique, and each $t\ge0$ falls into exactly one of the two categories: 
\begin{enumerate}[(i)]
\item For some $\delta>0$ and a strictly increasing function $\psi:\,(0,\delta]\to(0,\infty)$, the density of $X_{t-}\,\mathbf{1}_{\{\tau\ge t\}}$ on $(0,\delta]$ is bounded above by $\frac{1}{\alpha}-\psi$. In particular, $\Lambda_t=\Lambda_{t-}$.  
\item For some $\delta>0$ and a non-decreasing function $\psi:(0,\delta]\to[0,\infty)$, the density of $X_{t-}\,\mathbf{1}_{\{\tau\ge t\}}$ on $(0,\delta]$ is bounded below by $\frac{1}{\alpha}+\psi$. In particular, $\Lambda_t>\Lambda_{t-}$, and the density of $X_t\,\mathbf{1}_{\{\tau>t\}}=(X_{t-}-(\Lambda_t-\Lambda_{t-}))\,\mathbf{1}_{\{\tau>t\}}$ on $(0,\widetilde{\delta}]$ is bounded above by $\frac{1}{\alpha}-\widetilde{\psi}$, for some $\widetilde{\delta}>0$ and a strictly increasing function $\widetilde{\psi}:\,(0,\widetilde{\delta}]\to(0,\infty)$. 
\end{enumerate}
 \end{prop}

We refer to $t\ge0$ of categories (i) and (ii) in Proposition \ref{thm1} simply as \textit{continuity} and \textit{discontinuity points}.  



\subsection{Regularity estimates in right time neighborhoods} \label{rightn}

In this subsection, we consider an arbitrary $t\ge0$ and provide explicit estimates on the densities of $X_s\,\mathbf{1}_{\{\tau>s\}}$, $s\in(t,t+\teps]$ and $\min_{0\le m\le n} X_{t+m\Delta}\,\mathbf{1}_{\{\tau>t\}}$, $n=0,1,\ldots,\lfloor\teps/\Delta\rfloor$ in the vicinity of $0$, for suitable (explicit) $\teps>0$. These estimates play a crucial role in our proof of Theorem \ref{thm:loc_conv}.  

\begin{prop}\label{prop:density_est}
Suppose Assumption \ref{ass}. Let $t\ge0$ be a continuity (discontinuity resp.) point, with the density of $X_{t-}\,\mathbf{1}_{\{\tau\ge t\}}$ ($X_t\,\mathbf{1}_{\{\tau>t\}}$ resp.) on an interval $(0,\delta]$ being bounded above by $\frac{1}{\alpha}-\psi$ for a strictly increasing function $\psi:\,(0,\delta]\to(0,\infty)$. Then, with $\Psi(x):=\int_0^x \psi(y)\,\mathrm{d}y$, $x\in(0,\delta]$:
\begin{enumerate}[(a)]
	\item The density of $X_s\,\mathbf{1}_{\{\tau>s\}}$ on $(0,\delta/6)$ cannot exceed $\frac{1}{\alpha}-\frac{\psi}{2}$ for any 
	\begin{equation}
	t<s<t+\frac{\pi}{2\|f\|_\infty^2}\,\Psi(\delta/6)^2.
	\end{equation} 
	\item The density of $\min_{0\le m\le n} X_{t+m\Delta}\,\mathbf{1}_{\{\tau>t\}}$ on $(-\delta,0)$ is bounded above by
    \begin{equation}\label{<0UBD}
    \begin{cases}
    \frac{1}{\alpha}-c_1 & \text{if}\;\;\;\psi(0):=\lim_{y\downarrow 0} \psi(y)>0, \\
    \frac{1}{\alpha}-c_2\,\Phi\big(q-\frac{c_3}{\Psi(-x)}\big) & \text{if}\;\;\;\psi(0)=0,
    \end{cases}
    \end{equation}	
    $n=1,2,\ldots,\lfloor\teps/\Delta\rfloor$, $\Delta>0$, where $\Phi$, $c_1$, $c_2$, $c_3$, $q$ and $\teps$ are as in Theorem \ref{thm:loc_conv}. 
\end{enumerate}
\end{prop}

Our~proof~of~Proposition~\ref{prop:density_est}~relies~on~an~(explicit)~\textit{a}~\textit{priori}~estimate~on~the~increments~of~$\Lambda$. 

\begin{lem}\label{lemma:a priori}
In the situation of Proposition \ref{prop:density_est}, we have
\begin{equation}
\Lambda_s-\Lambda_t \le \Psi^{-1}\big(\|f\|_\infty\,\sqrt{2/\pi}\,\sqrt{s-t}\big),
\end{equation}
for all 
\begin{equation}\label{s_def}
t<s<\min\bigg(t+\frac{\pi}{2\|f\|_\infty^2}\,\Psi(\delta)^2,\,\inf\{r\!>\!t:\Lambda_r\!>\!\Lambda_{r-}\}\!\bigg).
\end{equation} 
\end{lem}

\noindent\textbf{Proof.} We use the definition of $\Lambda$ in \eqref{eq::init2}, its monotonicity, and the bounds $\frac{1}{\alpha}-\psi$, $\|f\|_\infty$ on the density of $X_t\,\mathbf{1}_{\{\tau>t\}}$ on $(0,\delta]$, $(0,\infty)$, respectively, to obtain
\begin{equation*}
\begin{split}
&\;\Lambda_s-\Lambda_t = \alpha\pp\Big(\tau>t,\,\inf_{r\in(t,s]} \big(X_t+(B_r-B_t)-(\Lambda_r-\Lambda_t)\big) \le 0\Big) \\
& \le \alpha\pp\Big(\tau>t,\,\inf_{r\in(t,s]} \big(X_t+(B_r-B_t)-(\Lambda_s-\Lambda_t)\big) \le 0\Big) \\
& \le \alpha\int_0^{\Lambda_s-\Lambda_t} \left(\frac{1}{\alpha}-\psi(x)\right)\,\mathrm{d}x+\alpha\int_{\Lambda_s-\Lambda_t}^\infty \|f\|_\infty\,\pp\Big(x+\inf_{r\in(t,s]} (B_r-B_t)-(\Lambda_s-\Lambda_t)\le 0\Big)\,\mathrm{d}x \\
&=\Lambda_s-\Lambda_t - \alpha\Psi(\Lambda_s-\Lambda_t)
+\alpha\|f\|_\infty\,\sqrt{2/\pi}\,\sqrt{s-t},
\end{split}
\end{equation*}
as long as $\Lambda_s-\Lambda_t<\delta$. In other words, $\Lambda_s-\Lambda_t<\delta$ implies
\begin{equation}
\Lambda_s-\Lambda_t\le\Psi^{-1}\big(\|f\|_\infty\,\sqrt{2/\pi}\,\sqrt{s-t}\big). 
\end{equation}
The lemma readily follows upon noting that \eqref{s_def} enforces $\Psi^{-1}\big(\|f\|_\infty\,\sqrt{2/\pi}\,\sqrt{s-t}\big)<\delta$. \qed

\medskip

We are now ready to give the proof of Proposition \ref{prop:density_est}. 

\medskip

\noindent\textbf{Proof of Proposition \ref{prop:density_est}. (a).} The density of $X_s\,\mathbf{1}_{\{\tau>s\}}$ on $(0,\infty)$ is bounded above by that of $\big(X_t+(B_s-B_t)-(\Lambda_s-\Lambda_t)\big)\,\mathbf{1}_{\{\tau>t\}}$. Recalling the notation $p(t,\cdot)$ for the density of $X_t\,\mathbf{1}_{\{\tau>t\}}$ on $(0,\infty)$ and writing $\varphi_{s-t}$ for the Gaussian density of mean $0$ and variance $s-t$, we can estimate the density of $\big(X_t+(B_s-B_t)-(\Lambda_s-\Lambda_t)\big)\,\mathbf{1}_{\{\tau>t\}}$ on $(0,\delta/6)$ by
\begin{equation*}
\begin{split}
&\,\int_{-\infty}^{-\frac{2}{3}\delta} p(t,x-y+\Lambda_s-\Lambda_t)\,\varphi_{s-t}(y)\,\mathrm{d}y
+\int_{-\frac{2}{3}\delta}^{x+\Lambda_s-\Lambda_t} p(t,x-y+\Lambda_s-\Lambda_t)\,\varphi_{s-t}(y)\,\mathrm{d}y \\
&\le\,\|f\|_\infty\,\Phi\bigg(\!-\frac{2\delta}{3\sqrt{s-t}}\bigg)
+\int_{-\frac{2}{3}\delta}^0 \bigg(\frac{1}{\alpha}-\psi(x-y+\Lambda_s-\Lambda_t)\!\bigg)\,\varphi_{s-t}(y)\,\mathrm{d}y \\
&\quad\;+\int_0^{x+\Lambda_s-\Lambda_t} \bigg(\frac{1}{\alpha}-\psi(x-y+\Lambda_s-\Lambda_t)\!\bigg)\,\varphi_{s-t}(y)\,\mathrm{d}y \\
&\le \,\|f\|_\infty\,\Phi\bigg(\!-\frac{2\delta}{3\sqrt{s-t}}\bigg)+\frac{1}{2}\,\bigg(\frac{1}{\alpha}-\psi(x)\!\bigg)
+\frac{1}{\alpha}\,\bigg(\frac{1}{2}-\Phi\bigg(\!-\frac{x+\Lambda_s-\Lambda_t}{\sqrt{s-t}}\bigg)\!\bigg),
\end{split}
\end{equation*}
provided that $\Lambda_s-\Lambda_t\le\delta/6$ (since then $x-y+\Lambda_s-\Lambda_t\in(0,\delta]$ for all $-\frac{2}{3}\delta\le y<x+\Lambda_s-\Lambda_t$) and where $\Phi$ is the standard Gaussian cumulative distribution function. Using $x\in(0,\delta/6)$ and the assumption $\Lambda_s-\Lambda_t\le\delta/6$ again we obtain the further upper bound 
\begin{equation*}
\frac{1}{\alpha}-\frac{\psi(x)}{2}
+\,\|f\|_\infty\,\Phi\bigg(\!-\frac{2\delta}{3\sqrt{s-t}}\bigg)-\frac{1}{\alpha}\,\Phi\bigg(\!-\frac{\delta}{3\sqrt{s-t}}\bigg).
\end{equation*} 
In view of the elementary inequality $\Phi(-2y)\le e^{-3y^2/2}\,\Phi(-y)$, $y>0$, it holds 
\begin{equation*}
\begin{split}
\frac{1}{\alpha}-\frac{\psi(x)}{2}
+\,\|f\|_\infty\,\Phi\bigg(\!-\frac{2\delta}{3\sqrt{s-t}}\bigg)-\frac{1}{\alpha}\,\Phi\bigg(\!-\frac{\delta}{3\sqrt{s-t}}\bigg)\le \frac{1}{\alpha}-\frac{\psi(x)}{2},\\ 
t<s<t+\frac{\delta^2}{6\log\max(\alpha\|f\|_\infty,1)}.
\end{split}
\end{equation*} 

\smallskip

Since $\Lambda_s-\Lambda_t\le\delta/6$ for all 
\begin{equation*}
t<s<\min\bigg(t+\frac{\pi}{2\|f\|_\infty^2}\,\Psi(\delta/6)^2,\,\inf\{r\!>\!t:\Lambda_r\!>\!\Lambda_{r-}\}\!\bigg)
\end{equation*}
by Lemma \ref{lemma:a priori}, we obtain the desired density estimate for all 
\begin{equation}\label{min3terms}
t<s<\min\bigg(t+\frac{\pi}{2\|f\|_\infty^2}\,\Psi(\delta/6)^2,\,\inf\{r\!>\!t:\Lambda_r\!>\!\Lambda_{r-}\},\,t+\frac{\delta^2}{6\log\max(\alpha\|f\|_\infty,1)}\bigg).
\end{equation}
In view of $\psi\le\frac{1}{\alpha}$, one has $\Psi(\delta/6)\le\delta/(6\alpha)$, and we conclude by elementary calculations that the first term in the minimum is less or equal to the third term. It remains to show that $\inf\{r\!>\!t:\Lambda_r\!>\!\Lambda_{r-}\}$ cannot attain the minimum in \eqref{min3terms}. Assuming the opposite and letting $\widetilde{t}:=\inf\{r\!>\!t:\Lambda_r\!>\!\Lambda_{r-}\}<\infty$ we infer that the density of $X_{\widetilde{t}-}\,\mathbf{1}_{\{\tau\ge\widetilde{t}\}}=\lim_{r\uparrow\widetilde{t}} X_r\,\mathbf{1}_{\{\tau>r\}}$ on $(0,\delta/6)$ is bounded above by $\frac{1}{\alpha}-\frac{\psi}{2}$, thus $\Lambda_{\widetilde{t}}=\Lambda_{\widetilde{t}-}$ due to the minimality property \eqref{physical}. This and the right-continuity of $\Lambda$ at $\widetilde{t}$ contradict $\widetilde{t}=\inf\{r\!>\!t:\Lambda_r\!>\!\Lambda_{r-}\}<\infty$. 

\medskip

\noindent\textbf{(b).} Fix $\Delta>0$, $1\le n\le\lfloor\teps/\Delta\rfloor$ and denote the law of $\min_{0\le m\le n} (B_{t+m\Delta}-B_t-\Lambda_{t+m\Delta}+\Lambda_t)$ by $\mu$. The density of $\min_{0\le m\le n} X_{t+m\Delta}\,\mathbf{1}_{\{\tau>t\}}$ on $(-\delta,0)$ is $\int_{(-\infty,x)} p(t,x-y)\,\mu(\mathrm{d}y)$, which can be bounded above by 
\begin{equation}\label{quantity}
\int_{(x-\delta,x)} \left(\frac{1}{\alpha}-\psi(x-y)\right)\,\mu(\mathrm{d}y) + \|f\|_\infty\,\mu\big((-\infty,x-\delta]\big).
\end{equation}
Next, we observe the inequalities 
\begin{equation}\label{ineq:comp}
\begin{split}
B_{t+n\Delta}-B_t-\Lambda_{t+n\Delta}+\Lambda_t &\ge\min_{0\le m\le n} (B_{t+m\Delta}-B_t-\Lambda_{t+m\Delta}+\Lambda_t) \\
&\ge \min_{t\le s\le t+n\Delta} (B_s-B_t) - (\Lambda_{t+n\Delta}-\Lambda_t)
\end{split}
\end{equation}
and introduce the notations $\widetilde{\varphi}$ and $\widetilde{\Phi}$ for the Gaussian density and cumulative distribution function of mean $-\Lambda_{t+n\Delta}+\Lambda_t$ and variance $n\Delta$. At this point, we distinguish two cases: $\mu\big((-\infty,x]\big)<\frac{1}{2}$ and $\mu\big((-\infty,x]\big)\ge\frac{1}{2}$. 

\medskip

In the first case, we apply $\mu\big((-\infty,x]\big)<\frac{1}{2}$, the second inequality in \eqref{ineq:comp} and the reflection principle for Brownian motion to estimate the density of $\min_{0\le m\le n} X_{t+m\Delta}\,\mathbf{1}_{\{\tau>t\}}$ on $(-\delta,0)$ further by 
\begin{equation*}
\frac{1}{2\alpha}+2\|f\|_\infty\,\widetilde{\Phi}(x-\delta)
=\frac{1}{2\alpha}+2\|f\|_\infty\,\Phi\bigg(\frac{x-\delta+\Lambda_{t+n\Delta}-\Lambda_t}{\sqrt{n\Delta}}\bigg).
\end{equation*}
Since $x\le0$ and $\Lambda_{t+n\Delta}-\Lambda_t\le\delta/6$ (see the last paragraph in the proof of part (a)), the latter estimate is less or equal to $\frac{1}{2\alpha}+2\|f\|_\infty\,\Phi\big(\!-\frac{5\delta}{6\sqrt{\teps}}\big)$. An elementary calculation shows that the range of $\teps>0$ for which this is less than $\frac{1}{\alpha}$ contains the interval specified by \eqref{teps:UBD}. 

\medskip

In the second case, we upper bound the quantity in \eqref{quantity} by
\begin{equation*}
\int_{\rr} \Big( \frac{1}{\alpha}-\psi(x-y)\Big)\,\mathbf{1}_{(x-\delta,x)}(y)+\frac{1}{\alpha}\,\mathbf{1}_{[x,\infty)}(y)\,\mu(\mathrm{d}y)+\|f\|_\infty\,\mu\big((-\infty,x-\delta]\big).
\end{equation*}
Using the inequalities in \eqref{ineq:comp} and the reflection principle for Brownian motion, we estimate the density of $\min_{0\le m\le n} X_{t+m\Delta}\,\mathbf{1}_{\{\tau>t\}}$ on $(-\delta,0)$ further by
\begin{equation}\label{expr}
\int_{\rr} \bigg(\!\Big( \frac{1}{\alpha}-\psi(x-y)\Big)\,\mathbf{1}_{(x-\delta,x)}(y)+\frac{1}{\alpha}\,\mathbf{1}_{[x,\infty)}(y)\!\bigg)
\,\widetilde{\varphi}(y)\,\mathrm{d}y +2\|f\|_\infty\,\widetilde{\Phi}(x-\delta).
\end{equation}
Next, we exploit that for any $z\in[0,\delta)$, it holds $\psi(x-y)\ge\psi(z)$, $y\in(x-\delta,x-z)$ and $\psi(x-y)\ge0$, $y\in[x-z,x)$. This allows to upper bound the expression in \eqref{expr} by
\begin{equation*}
\begin{split}
& \;\frac{1}{\alpha}\big(1-\widetilde{\Phi}(x-\delta)\big)-\psi(z)\big(\widetilde{\Phi}(x-z)-\widetilde{\Phi}(x-\delta)\big)+2\|f\|_\infty\,\widetilde{\Phi}(x-\delta) \\
& =\frac{1}{\alpha}-\psi(z)\,\widetilde{\Phi}(x-z)+\Big(2\|f\|_\infty+\psi(z)-\frac{1}{\alpha}\Big)\,\widetilde{\Phi}(x-\delta).
\end{split}
\end{equation*}
In view of the elementary inequality 
\begin{equation*}
\widetilde{\Phi}(x-\delta)\le\widetilde{\Phi}(x-z)\,\exp\bigg(\!-\frac{(\delta-z)(\delta+z-2x-2\Lambda_{t+n\Delta}+2\Lambda_t)}{2n\Delta}\bigg),
\end{equation*}
the latter cannot exceed 
\begin{equation}\label{case2bnd}
\frac{1}{\alpha}-\widetilde{\Phi}(x-z)\,\bigg(\psi(z)-\Big(2\|f\|_\infty+\psi(z)-\frac{1}{\alpha}\Big)_+\,\exp\bigg(\!-\frac{(\delta-z)(\delta+z-2x-2\Lambda_{t+n\Delta}+2\Lambda_t)}{2n\Delta}\bigg)\!\bigg).
\end{equation}

\smallskip

Combining $\mu\big((-\infty,x]\big)\ge\frac{1}{2}$, the second inequality in \eqref{ineq:comp} and the reflection principle for Brownian motion we find $\widetilde{\Phi}(x)\ge\frac{1}{4}$. In other words, $x\ge -\Lambda_{t+n\Delta}+\Lambda_t+\sqrt{n\Delta}\,q$, where $q$ is the $\frac{1}{4}$-quantile of the standard Gaussian distribution. Therefore, writing $\Phi$ for the standard Gaussian cumulative distribution function again, we conclude 
\begin{equation}\label{Gaussianbnd}
\widetilde{\Phi}(x-z)=\Phi\bigg(\frac{x-z+\Lambda_{t+n\Delta}-\Lambda_t}{\sqrt{n\Delta}}\bigg)
\ge\Phi\bigg(\frac{-z+\sqrt{n\Delta}\,q}{\sqrt{n\Delta}}\bigg)
=\Phi\bigg(\!-\frac{z}{\sqrt{n\Delta}}+q\bigg).
\end{equation}
In addition, recalling part (a), the minimality property \eqref{physical} and Lemma \ref{lemma:a priori} we obtain
\begin{equation*}
x\ge -\Lambda_{t+n\Delta}+\Lambda_t+q\,\sqrt{n\Delta}
\ge -\Psi^{-1}\big(\|f\|_\infty\,\sqrt{2/\pi}\,\sqrt{n\Delta}\big)+q\,\sqrt{n\Delta}.
\end{equation*}
This yields 
\begin{equation*}
\Psi^{-1}\big(\Psi(-x)\big)-\Psi^{-1}\big(\|f\|_\infty\,\sqrt{2/\pi}\,\sqrt{n\Delta}\big)\le -q\,\sqrt{n\Delta}
\end{equation*}
and, since $(\Psi^{-1})'\ge 1/\psi(\delta)$, 
\begin{equation}\label{SDbnd}
\Psi(-x)-\|f\|_\infty\,\sqrt{2/\pi}\,\sqrt{n\Delta}\le -q\,\sqrt{n\Delta}\,\psi(\delta)\quad\Longleftrightarrow\quad
\sqrt{n\Delta}\ge\frac{\Psi(-x)}{\|f\|_\infty\,\sqrt{2/\pi}-q\psi(\delta)}. 
\end{equation} 

\smallskip

Finally, we choose $z=0$ if $\psi(0):=\lim_{y\downarrow0} \psi(y)>0$ and $z=\frac{\delta}{2}$ if $\psi(0)=0$, and combine \eqref{case2bnd}, $\widetilde{\Phi}(x)\ge\frac{1}{4}$, \eqref{Gaussianbnd} and \eqref{SDbnd} to get the estimate
\begin{equation*}
\begin{cases}
\frac{1}{\alpha}-\frac{1}{4}\,\Big(\psi(0)-\big(2\|f\|_\infty+\psi(0)-\frac{1}{\alpha}\big)_+\,\exp\Big(\!-\frac{\delta(\delta-2x-2\Lambda_{t+n\Delta}+2\Lambda_t)}{2n\Delta}\Big)\Big) & \text{if}\;\;\;\psi(0)>0, \\
\frac{1}{\alpha}-\Phi\Big(q-\frac{\delta(\|f\|_\infty\,\sqrt{2/\pi}-q\psi(\delta))}{2\Psi(-x)}\Big) & \quad \\
\quad\;\;\;\,\cdot\Big(\psi(\delta/2)-\big(2\|f\|_\infty+\psi(\delta/2)-\frac{1}{\alpha}\big)_+\exp\Big(\!-\frac{\delta(3\delta/2-2x-2\Lambda_{t+n\Delta}+2\Lambda_t)}{4n\Delta}\Big)\Big) &
 \text{if}\;\;\;\psi(0)=0
\end{cases}
\end{equation*}
on the density of $\min_{0\le m\le n} X_{t+m\Delta}\,\mathbf{1}_{\{\tau>t\}}$ on $(-\delta,0)$ in the second case. (Note that the respective long parenthesis is non-negative thanks to $x\le 0$, $\Lambda_{t+n\Delta}-\Lambda_t\le\delta/6$, $n\Delta\le\teps$ and \eqref{teps:UBD}, as well as to $\psi(\delta/2)\le\frac{1}{\alpha}$ and $\Psi(\delta/6)\le\psi(\delta/2)\delta/6$ when $\psi(0)=0$.) We conclude the proof by using $x\le 0$, $\Lambda_{t+n\Delta}-\Lambda_t\le\delta/6$ and $n\Delta\le\teps$, and then putting the resulting bound together with the one found in the first case. (Hereby, $\frac{1}{4}\le\frac{1}{2}$ and $\psi\le\frac{1}{\alpha}$ are exploited.) \qed

\medskip

By combining Proposition \ref{prop:density_est}(a) with Lemma \ref{lemma:a priori} we can control the modulus of continuity of $\Lambda$ on right time neighborhoods. 

\begin{cor}\label{cor:mod_of_cont}
In the situation of Proposition \ref{prop:density_est}, it holds 
\begin{equation}
\Lambda_{s_2}-\Lambda_{s_1}\le 2\Psi^{-1}\big(\|f\|_\infty\,\sqrt{2/\pi}\,\sqrt{s_2-s_1}\big),\quad t\le s_1<s_2\le t+\frac{\pi}{8\|f\|_\infty^2}\,\Psi(\delta/6)^2.
\end{equation}
\end{cor}

\noindent\textbf{Proof.} It suffices to apply Lemma \ref{lemma:a priori} with $t=s_1$ upon recalling the density estimate of Proposition \ref{prop:density_est}(a) with $s=s_1$. \qed

\subsection{Useful definitions}\label{subsec:M1}

In this subsection, we prepare the definitions of the Skorokhod M1 convergence for non-decreasing functions and of some auxiliary time-stepping schemes used in the proof of Theorem \ref{thm::main}. 

%

\begin{dfn}[see, e.g., \cite{Whi}, Corollary 12.5.1]\label{M1def}
Let $T>0$. A sequence $(L_k)_{k\in\mathbb{N}}$ of non-decreasing functions in $D([0,T])$ is said to \emph{converge in M1 sense} to some $L\in D([0,T])$ if and only if $\lim_{k\to\infty} L_k(t)=L(t)$ for all $t$ in a dense subset of $[0,T]$ that includes $0$ and $T$. 
\end{dfn}

\begin{dfn}\label{defn::discrete_general}
\begin{enumerate}[(a)]	
\item Given $t\ge0$, $\Delta>0$, a physical solution $\Lambda$ and the associated $X$, we define $X^{t;\Delta}_{t+n\Delta}$, $n=0,\,1,\,\ldots$ and $\Lambda^{t;\Delta}_{t+n\Delta}$, $n=1,\,2,\,\ldots$ recursively by $X^{t;\Delta}_t=X_t$, 
\begin{equation}
\begin{split}
& \Lambda^{t;\Delta}_{t+n\Delta}=\alpha\pp\Big(\inf_{0\le s\le t} X_s\le 0\;\,\text{or}\,
\min_{0\le m\le n-1} X^{t;\Delta}_{t+m\Delta}<0\Big),\;n=1,2,\ldots\;\;\text{and} \\ 
& X^{t;\Delta}_{t+n\Delta}=X_{0-}+B_{t+n\Delta}-\Lambda^{t;\Delta}_{t+n\Delta},\;n=1,2,\ldots.  
\end{split}
\end{equation}
We extrapolate by setting $\Lambda^{t;\Delta}_{0-}=0$, $\Lambda^{t;\Delta}_s=\Lambda_s$, $s\in[0,t]$ and $\Lambda^{t;\Delta}_s=\Lambda^{t;\Delta}_{\max\{t+n\Delta:\,t+n\Delta\le s\}}$, $s>t$, as well as $X^{t;\Delta}_{0-}=X_{0-}$ and $X^{t;\Delta}_s = X_{0-} + B_s - \Lambda^{t;\Delta}_s$, $s\ge 0$. 
\item Given a sequence $\CTT$ of the form $0=t_0<t_1<t_2<\cdots$, we let $X^{\CTT}_{t_n}$, $n=0,\,1,\,\ldots$ and $\Lambda^{\CTT}_{t_n}$, $n=1,\,2,\,\ldots$ be defined recursively by $X^{\CTT}_{t_0}=X_0$, 
\begin{equation}
\Lambda^{\CTT}_{t_n}=\alpha\pp\Big(\min_{0\le m\le n-1} X^{\CTT}_{t_m}<0\Big),\;n=1,2,\ldots\;\;\text{and}\;\;\, X^{\CTT}_{t_n}=X_{0-}+B_{t_n}-\Lambda^{\CTT}_{t_n},\;n=1,2,\ldots.  
\end{equation}
We extrapolate by setting $\Lambda^{\CTT}_{0-}=0$, $\Lambda^{\CTT}_0=\alpha\pp(X_0\le0)$ and $\Lambda^{\CTT}_t=\Lambda^{\CTT}_{\max\{t_n:\,t_n\le t\}}$, $t>0$, as well as $X^{\CTT}_{0-}=X_{0-}$ and $X^{\CTT}_t = X_{0-} + B_t - \Lambda^{\CTT}_t$, $t\ge0$. 
\end{enumerate}
\end{dfn}

\begin{rmk}\label{rmk::discrete_general_prelim}
If $\CTT,\widetilde{\CTT}$ are sequences as in Definition \ref{defn::discrete_general}(b) and $\CTT$ is a subsequence of $\widetilde{\CTT}$, then $\La^{\CTT}\le\La^{\widetilde{\CTT}}\le\La$, by induction. 
\end{rmk}


\section{Proof of Theorem \ref{thm:loc_conv}} \label{sec:loc_conv}

In this section, we establish the following generalization of Theorem \ref{thm:loc_conv} to arbitrary $t\ge0$ (instead of only $t=0$), which is also used in the proof of Theorem \ref{thm::main} below. 

\begin{prop}\label{prop:loc_conv}
In the situation of Proposition \ref{prop:density_est}, let $\frac{1}{\alpha}-\widetilde{\psi}:\,(-\delta,0)\to\big[0,\frac{1}{\alpha}\big)$ be the function in \eqref{<0UBD} and $\widetilde{\Psi}(x):=\int_0^x \widetilde{\psi}(-y)\,\mathrm{d}y$, $x\in(0,\delta)$. Then, it holds 
\begin{equation}\label{loc_rate_conv}
\begin{split}
\sup_{s\in[t,t+\teps]} |\Lambda^{t;\Delta}_s-\Lambda_s| \le \widetilde{\Psi}^{-1}\Big(\|f\|_\infty\Big(192\sqrt{2\Delta\log(\lceil\teps/\Delta\rceil/2)} 
+2\Psi^{-1}\big(2\sqrt{\pi}^{-1}\|f\|_\infty\sqrt{\Delta}\big)\!\Big)\!\Big)\quad\;\;\; \\
+\max(\alpha\|f\|_\infty,1)\Big(192\sqrt{2\Delta\log(\lceil\teps/\Delta\rceil/2)} 
\!+\!2\Psi^{-1}\big(2\sqrt{\pi}^{-1}\|f\|_\infty\sqrt{\Delta}\big)\!\Big)
\end{split} 
\end{equation}
for all $\teps\!>\!0$ adhering to \eqref{teps:UBD} and all $\Delta>0$ small enough.
\end{prop}

\begin{rmk}
For a discontinuity point $t>0$, the function $\psi$ in Proposition \ref{prop:density_est} can be chosen as $cx^a$ for some $a\in{\mathbb N}$ and $c>0$ (see, e.g., \cite[proof of Theorem 3.3]{LS2}). Then, the argument of $\widetilde{\Psi}^{-1}$ and the term on the second line in \eqref{loc_rate_conv} have the order $1/(2(a+1))$ in $\Delta$. As touched upon in Remark \ref{rm_order}, we conjecture this to reflect the true order of the left-hand side in~\eqref{loc_rate_conv}, based on the numerical simulations described in Subsection \ref{sec5.3}.
\end{rmk}

\noindent\textbf{Proof of Proposition \ref{prop:loc_conv}.} Clearly, it suffices to check that each of $\Lambda_{t+n\Delta}-\Lambda^{t;\Delta}_{(t+n\Delta)-}$, $n=1,2,\ldots,\lceil\teps/\Delta\rceil$ does not exceed the right-hand side of \eqref{loc_rate_conv}. To this end, we argue by induction over $n$. For $n=1$, the result follows from $\Lambda^{t;\Delta}_{(t+\Delta)-}=\Lambda_t$ and Corollary \ref{cor:mod_of_cont}, so we turn to $n\ge 2$. We then have
\begin{equation}\label{eq:2diff}
\begin{split}
\Lambda_{t+n\Delta}\!-\!\Lambda^{t;\Delta}_{(t+n\Delta)-}&=\alpha\pp\Big(\inf_{0\le s\le t+n\Delta} X_s\le 0\Big)
\!-\!\alpha\pp\Big(\inf_{0\le s\le t} X_s\le 0\;\,\text{or}\;\min_{0\le m\le n-2} X^{t;\Delta}_{t+m\Delta}<0\Big) \\
&=\alpha\pp\Big(\tau>t,\,\inf_{t<s\le t+n\Delta} X_s\le 0\Big)
-\alpha\pp\Big(\tau>t,\,\min_{0\le m\le n-2} X^{t;\Delta}_{t+m\Delta}<0\Big) \\
& =\alpha\pp\Big(\tau>t,\,\inf_{t<s\le t+n\Delta} X_s\le 0\Big)
-\alpha\pp\Big(\tau>t,\,\min_{0\le m\le n-2} X_{t+m\Delta}<0\Big) \\
&\quad+\alpha\pp\Big(\tau\!>\!t,\,\min_{0\le m\le n-2} X_{t+m\Delta}<0\Big)
\!-\!\alpha\pp\Big(\tau\!>\!t,\,\min_{0\le m\le n-2} X^{t;\Delta}_{t+m\Delta}<0\Big).
\end{split}
\end{equation}

\smallskip

The difference on the third line in \eqref{eq:2diff} is bounded above by
\begin{equation*}
\begin{split}
& \;\alpha\pp\Big(0<X_t\,\mathbf{1}_{\{\tau>t\}}+\min_{0\le m\le n-2} (B_{t+m\Delta}-B_t-\Lambda_{t+m\Delta}+\Lambda_t) \\
&\qquad\;\;\,\le\sup_{{t\le s_1<s_2\le t+n\Delta}\atop{s_2-s_1\le 2\Delta}} |B_{s_2}-B_{s_1}|+\sup_{{t\le s_1<s_2\le t+n\Delta}\atop{s_2-s_1\le 2\Delta}} (\Lambda_{s_2}-\Lambda_{s_1})\!\Big) \\
& \le \alpha\|f\|_\infty\bigg(\E\Big[\sup_{{t\le s_1<s_2\le t+n\Delta}\atop{s_2-s_1\le 2\Delta}} |B_{s_2}-B_{s_1}|\Big]+\sup_{{t\le s_1<s_2\le t+n\Delta}\atop{s_2-s_1\le 2\Delta}} (\Lambda_{s_2}-\Lambda_{s_1})\!\bigg) \\
& \le \alpha\|f\|_\infty\Big(192\sqrt{2\Delta\log(\lceil\teps/\Delta\rceil/2)} 
+2\Psi^{-1}\big(2\sqrt{\pi}^{-1}\|f\|_\infty\sqrt{\Delta}\big)\!\Big), 
\end{split}
\end{equation*}
where we have used  \cite[Lemma 4 and Remark 3]{Fi} and Corollary \ref{cor:mod_of_cont}. The difference on the fourth line in \eqref{eq:2diff} can be estimated by 
\begin{equation*}
\begin{split}
&\;\alpha\pp\Big(-\max_{0\le m\le n-2}\,\big(\Lambda_{t+m\Delta}-\Lambda^{t;\Delta}_{(t+m\Delta)-}\big)\le\min_{0\le m\le n-2} X_{t+m\Delta}\,\mathbf{1}_{\{\tau>t\}}<0\Big) \\
&\le \max_{0\le m\le n-2}\,\big(\Lambda_{t+m\Delta}-\Lambda^{t;\Delta}_{(t+m\Delta)-}\big)
-\alpha\widetilde{\Psi}\Big(\max_{0\le m\le n-2}\big(\Lambda_{t+m\Delta}-\Lambda^{t;\Delta}_{(t+m\Delta)-}\big)\!\Big)
\end{split}
\end{equation*}
thanks to the induction hypothesis and Proposition \ref{prop:density_est}(b). 

\medskip

All in all, we have obtained
\begin{equation*}
\begin{split}
\Lambda_{t+n\Delta}\!-\!\Lambda^{t;\Delta}_{(t+n\Delta)-}\le  & \,\max_{0\le m\le n-2} \big(\Lambda_{t+m\Delta}\!-\!\Lambda^{t;\Delta}_{(t+m\Delta)-}\big)
-\alpha\widetilde{\Psi}\Big(\max_{0\le m\le n-2}\big(\Lambda_{t+m\Delta}\!-\!\Lambda^{t;\Delta}_{(t+m\Delta)-}\big)\!\Big) \\
&+\alpha\|f\|_\infty\Big(192\sqrt{2\Delta\log(\lceil\teps/\Delta\rceil/2)} 
+2\Psi^{-1}\big(2\sqrt{\pi}^{-1}\|f\|_\infty\sqrt{\Delta}\big)\!\Big).
\end{split}
\end{equation*}
We conclude by distinguishing the cases of how $\max_{0\le m\le n-2}\big(\Lambda_{t+m\Delta}-\Lambda^{t;\Delta}_{(t+m\Delta)-}\big)$ compares to the $\widetilde{\Psi}^{-1}$-term in \eqref{loc_rate_conv} and using the induction hypothesis. \qed  

\begin{rmk}\label{alt_ass}
It is worth noting that the proof of Proposition \ref{prop:loc_conv} (and thus, of Theorem~\ref{thm:loc_conv}) remains intact if instead of Assumption \ref{ass}(a), the boundedness of $f$ on $[0,\infty)$ together with the conclusions of Lemma \ref{lemma:a priori} and Corollary~\ref{cor:mod_of_cont} are assumed and $\psi:=\Psi':\,(0,\delta]\to\big(0,\frac{1}{\alpha}\big]$ is a strictly increasing function. Indeed, the proofs of Propositions \ref{prop:loc_conv} and \ref{prop:density_est}(b) can be then repeated word by word. 
\end{rmk}

\newcommand{\ep}{\hfill \ensuremath{\Box}}

\section{Proof of Theorem \ref{thm::main}} \label{sec:mainproof}

For the proof of the convergence globally in time (Theorem \ref{thm::main}), we prepare auxiliary error propagation bounds in Subsection \ref{sec::stab} and then show respectively that the convergence can be extended beyond continuity and discontinuity points of $\Lambda$ in Subsections \ref{sec::no_jump} and \ref{sec::special_case}.
%
\subsection{Auxiliary error propagation bounds} \label{sec::stab}
The following lemma gives an estimate concerning the continuity of the numerical solution with respect to the initial condition.
The statement is weaker than typical notions of stability as the ``constant'' in it may depend on $\Delta$, but this is sufficient for our purposes.

\begin{lem}\label{lem::discrete_compare}
For some $\Delta>0$ and $K\ge0$, let $(X^\Delta, \Lambda^\Delta)$ be given by Definition \ref{defn::discrete}, and let $\widetilde{X}^\Delta_{n\Delta}$, $n=0,\,1,\,\ldots$ and $\widetilde{\Lambda}^\Delta_{n\Delta}$, $n=0,\,1,\,\ldots$ be defined by $\widetilde{X}^\Delta_0=X^\Delta_0+K$, $\widetilde{\Lambda}^\Delta_0=\alpha\mathbb{P}(\widetilde{X}^\Delta_0\le0)$, 
\begin{equation}
\widetilde{\Lambda}^{\Delta}_{n\Delta}=\alpha\pp\Big(\min_{0\le m\le n-1} \widetilde{X}^\Delta_{m\Delta}<0\Big),\;n\ge 1\;\;\text{and}\;\; \widetilde{X}^\Delta_{n\Delta}=\widetilde{X}^\Delta_0+B_{n\Delta}-\widetilde{\Lambda}^\Delta_{n\Delta}
+\widetilde{\Lambda}^\Delta_0
,\;n\ge 1.  
\end{equation}
We extrapolate to $t>0$ by setting $\widetilde{\Lambda}^{\Delta}_t=\widetilde{\Lambda}^{\Delta}_{\lfloor t/\Delta \rfloor \Delta}$ and $\widetilde{X}^{\Delta}_t=\widetilde{X}^{\Delta}_{\lfloor t/\Delta \rfloor \Delta}$. If $X_{0-}\ge0$ admits a density on $[0,\infty)$, then for any $T>0$, there is some $C=C_\alpha(\Delta,T)\ge0$ non-decreasing in $\alpha$ and non-increasing in $\Delta$ such that
\begin{equation}
(\Lambda^\Delta_t-\Lambda^\Delta_0)-(\widetilde{\Lambda}^\Delta_t-\widetilde{\Lambda}^\Delta_0)\leq CK,\quad t\in[0,T].
\end{equation}
\end{lem}

\noindent\textbf{Proof.} For $n\!\ge\!1$, by $\min_{0\le m \le n-1} \widetilde{X}^\Delta_{m\Delta}\!\ge0\Leftrightarrow\widetilde{X}^\Delta_{m\Delta}\ge0,\,0\!\le\! m\!\le\! n\!-\!1$ and the union bound
\begin{equation*}
\begin{split}
\BPP \Big( \min_{0\le m \le n-1} X^\Delta_{m\Delta} < 0 \le \min_{0\le m \le n-1} \widetilde{X}^\Delta_{m\Delta} \Big) &\le\,
\BPP \Big(\exists\,0\le m\le n-1:\; X^\Delta_{m\Delta} < 0 \le \widetilde{X}^\Delta_{m\Delta}\Big) \\
&\le\,\sum_{m=0}^{n-1} \BPP \big( X^\Delta_{m\Delta} < 0 \le \widetilde{X}^\Delta_{m\Delta} \big).
\end{split}
\end{equation*}
Now, we set $K_n=(\Lambda^\Delta_{n\Delta}-\Lambda^\Delta_0)-(\widetilde{\Lambda}^\Delta_{n\Delta}-\widetilde{\Lambda}^\Delta_0)$, $n=1,\,2,\,\ldots\,$. Then, using the latter estimate and bounding the standard normal density by its maximum we infer
\begin{equation*}
\begin{split}
& \; K_n \le \alpha\pp\Big( \min_{0\le m \le n-1} X^\Delta_{m\Delta} < 0 \le \min_{0\le m \le n-1} \widetilde{X}^\Delta_{m\Delta} \Big)-\alpha\pp\big(X^\Delta_0\le 0<\widetilde{X}^\Delta_0\big) \\
& \le \sum_{m=1}^{n-1} \alpha\pp\big(X^\Delta_{m\Delta} < 0 \le \widetilde{X}^\Delta_{m\Delta} \big) 
=\sum_{m=1}^{n-1} \alpha\pp\Big( B_{m\Delta} \in \big[\widetilde{\Lambda}^\Delta_{m\Delta} -\widetilde{\Lambda}^\Delta_0-\widetilde{X}^\Delta_0,\Lambda^\Delta_{m\Delta} -\Lambda^\Delta_0-X^\Delta_0\big)\Big) \\
& \le \sum_{m=1}^{n-1} \frac{\alpha}{\sqrt{2\pi m\Delta}} (K_m + K).
\end{split}
\end{equation*}
By induction, we conclude that $K_n\le \widehat{C}_\alpha(\Delta,n)K$, $n=1,\,2,\, \ldots\,$, for suitable constants $\widehat{C}_\alpha(\Delta,n)\ge0$, $n=1,\,2,\,\ldots\,$.
It remains to take $C_\alpha(\Delta,T)=\max_{1\le n \le \lceil T/\Delta\rceil} \widehat{C}_\alpha(\Delta,n)$. \ep
	
\medskip

%

The extension proofs in Subsections \ref{sec::no_jump}, \ref{sec::special_case} rely on the following corollary of Lemma \ref{lem::discrete_compare}.

\begin{cor}\label{cor::discrete_compare}
For every $\Delta>0$ and $0\le r\le s$, there is some $C=C_\alpha(\Delta,s-r)\ge0$ non-increasing in $\Delta$ such that for any $N\ge0$ and sequence $\CTT$ of the form 
\begin{equation*}
0=t_0<\cdots<t_N=r<r+\Delta<r+2\Delta<\cdots, 
\end{equation*}
the function $\Lambda^{r;\Delta}$ from Definition \ref{defn::discrete_general}(a) and the function $\Lambda^{\CTT}$ from Definition 
\ref{defn::discrete_general}(b) satisfy
\begin{equation}
\Lambda^{r;\Delta}_{r+t}-\Lambda^{\CTT}_{r+t}\le C (\Lambda^{r;\Delta}_r-\Lambda^{\CTT}_r) ,\quad t\in[0,s-r],
\end{equation}
provided $X_{0-}\ge0$ admits a density on $[0,\infty)$. 
\end{cor}

\noindent\textbf{Proof.} We first note the comparison of increments
\begin{equation}\label{incr_comp}
\Lambda^{\CTT}_{r+t} - \Lambda^{\CTT}_r \ge \widetilde{\Lambda}^{r;\Delta}_{r+t} - \widetilde{\Lambda}^{r;\Delta}_r,\quad t\in[0,s-r], 
\end{equation} 
where $(\widetilde{X}^{r;\Delta}_{r+t},\widetilde{\Lambda}^{r;\Delta}_{r+t})$, $t\in[0,s-r]$ is defined by $\widetilde{X}^{r;\Delta}_r=X^{\CTT}_r$, $\widetilde{\Lambda}^{r;\Delta}_r=\alpha\pp\big(\inf_{0\le q\le r} X_q\le 0\big)$,
\begin{eqnarray*}
&& \widetilde{\Lambda}^{r;\Delta}_{r+n\Delta} = \alpha\pp\Big(\inf_{0\le q\le r} X_q\le 0\;\,\mathrm{or}\;
\min_{0\le m\le n-1} \widetilde{X}^{r;\Delta}_{r+m\Delta}<0\Big), \quad n=1,\,2,\,\ldots,\big\lfloor(s-r)/\Delta\big\rfloor,  \\
&& \widetilde{X}^{r;\Delta}_{r+n\Delta} = \widetilde{X}^{r;\Delta}_r + B_{r+n\Delta} - B_r - \widetilde{\Lambda}^{r;\Delta}_{r+n\Delta}+\widetilde{\Lambda}^{r;\Delta}_r,\quad n=1,\,2,\,\ldots,\,\big\lfloor(s-r)/\Delta\big\rfloor,
\end{eqnarray*}
and the extrapolations $\widetilde{\Lambda}^{r;\Delta}_{r+t} = \widetilde{\Lambda}^{r;\Delta}_{r+\lfloor t /\Delta\rfloor\Delta}$ and $\widetilde{X}^{r;\Delta}_{r+t} = \widetilde{X}^{r;\Delta}_{r+\lfloor t /\Delta\rfloor\Delta}$ for $t\in(0,s-r]$. Indeed, \eqref{incr_comp} is trivial for $t\in[0,\Delta)$ and holds for $t\in[\Delta,s-r]$, since 
\begin{equation}
\pp\Big(\min_{0\le m\le (N-1)_+} X^{\CTT}_{t_m}\!\ge 0,\,\min_{0\le m\le n-1} X^{\CTT}_{r+m\Delta}\!<0\Big) 
\ge \pp\Big(\inf_{0\le q\le r} X_q\!>0,\,\min_{0\le m\le n-1} \widetilde{X}^{r;\Delta}_{r+m\Delta}\!<0\Big),  
\end{equation}
$n:=\lfloor t /\Delta\rfloor=1,\,2,\,\ldots,\,\lfloor(s-r)/\Delta\rfloor$, as can be readily seen by induction.

\medskip

Conditioning on the event $\{\inf_{0\le q\le r} X_q>0\}$, which has probability $1-\Lambda_r/\alpha$, we get
\begin{eqnarray*}
&& \Lambda^{r;\Delta}_{r+t}-\Lambda^{r;\Delta}_r=(\alpha-\Lambda_r)\,\pp\Big(\min_{0\le m\le \lfloor t /\Delta\rfloor-1} X^{r;\Delta}_{r+m\Delta}\le 0\,\Big|\,\inf_{0\le q\le r} X_q>0\Big),\quad t\in[0,s-r], \\
&& \widetilde{\Lambda}^{r;\Delta}_{r+t}-\widetilde{\Lambda}^{r;\Delta}_r=(\alpha-\Lambda_r)\,\pp\Big(\min_{0\le m\le \lfloor t /\Delta\rfloor-1} \widetilde{X}^{r;\Delta}_{r+m\Delta}\le 0\,\Big|\,\inf_{0\le q\le r} X_q>0\Big),\quad t\in[0,s-r]. 
\end{eqnarray*} 
Thus, under the conditional measure $\pp(\,\cdot\,|\,\inf_{0\le q\le r} X_q>0)$, the pairs $(X^{r;\Delta}_{r+\cdot},\Lambda^{r;\Delta}_{r+\cdot}-\Lambda^{r;\Delta}_r)$ and $(\widetilde{X}^{r;\Delta}_{r+\cdot},\widetilde{\Lambda}^{r;\Delta}_{r+\cdot}-\widetilde{\Lambda}^{r;\Delta}_r)$ fit into the framework of Lemma \ref{lem::discrete_compare}, which yields the estimate
\begin{equation} \label{bound_from_lemma}
(\Lambda^{r;\Delta}_{r+t}-\Lambda^{r;\Delta}_r)-(\widetilde{\Lambda}^{r;\Delta}_{r+t}-\widetilde{\Lambda}^{r;\Delta}_r)\le C_{\alpha-\Lambda_r}(\Delta,s-r) K\le C_\alpha(\Delta,s-r) K,\quad t\in[0,s-r],\;\;\;
\end{equation}
where $K:=\Lambda^{r;\Delta}_r-\Lambda^{\CTT}_r=X^{\CTT}_r-X^{r;\Delta}_r=\widetilde{X}^{r;\Delta}_r-X^{r;\Delta}_r$. Using this definition of $K$, then \eqref{incr_comp}, and finally \eqref{bound_from_lemma}, we end up with
\begin{equation*}
\begin{split}
\Lambda^{r;\Delta}_{r+t}-\Lambda^{\CTT}_{r+t}
&=(\Lambda^{r;\Delta}_{r+t}-\Lambda^{r;\Delta}_r)-(\Lambda^{\CTT}_{r+t}-\Lambda^{\CTT}_r)+K \\
&\le (\Lambda^{r;\Delta}_{r+t}-\Lambda^{r;\Delta}_r)-(\widetilde{\Lambda}^{r;\Delta}_{r+t} - \widetilde{\Lambda}^{r;\Delta}_r)+K 
\le\big(C_\alpha(\Delta,s-r)+1\big)K. \qquad\qquad\quad\;\;\;\ep
\end{split}
\end{equation*}
	

\subsection{Extension beyond continuity points} \label{sec::no_jump}

We proceed to the following extension result.

%


\begin{prop}\label{pro::no_jump}
Suppose Assumption \ref{ass}. Let $T\ge0$ be a continuity point of $\Lambda$ such that $\Lambda^\Delta|_{[0,t]}\underset{\Delta\downarrow0}{\overset{\emph{M1}}{\longrightarrow}}\Lambda|_{[0,t]}$ for any continuity point $t\in(0,T)$ of $\Lambda$. Then,  $\Lambda^\Delta|_{[0,T+\teps]}\underset{\Delta\downarrow0}{\overset{\emph{M1}}{\longrightarrow}}\Lambda|_{[0,T+\teps]}$ for all $\teps>0$ small enough.
\end{prop}

\noindent\textbf{Proof.} We first argue that $\Lambda^\Delta|_{[0,T]}\underset{\Delta\downarrow0}{\overset{\text{M1}}{\longrightarrow}}\Lambda|_{[0,T]}$ if $T>0$. Thanks to the definition of the M1 convergence for monotone functions (Definition \ref{M1def}) and in view of the assumed convergence $\Lambda^\Delta|_{[0,t]}\underset{\Delta\downarrow0}{\overset{\text{M1}}{\longrightarrow}}\Lambda|_{[0,t]}$ for all continuity points $t\in(0,T)$ of $\Lambda$, it suffices to verify that $\lim_{\Delta\downarrow0} \Lambda^\Delta_T = \Lambda_T$. Taking a sequence $T_N\uparrow T$ as $N\to\infty$ satisfying $\lim_{\Delta\downarrow0}  \Lambda^\Delta_{T_N} = \Lambda_{T_N}$, $N=1,\,2,\,\ldots$ we obtain
\begin{equation}
\Lambda_T \ge \limsup_{\Delta\downarrow0}  \Lambda^\Delta_T 
\ge \liminf_{\Delta\downarrow0}  \Lambda^\Delta_T
\ge \lim_{\Delta\downarrow0}  \Lambda^\Delta_{T_N}=\Lambda_{T_N},\quad N=1,\,2,\,\ldots.
\end{equation}   
Since $\Lambda_{T_N}\uparrow \Lambda_T$ as $N\to\infty$, it must hold
\begin{equation}
\label{conv_T}
\Lambda_T=\lim_{\Delta\downarrow0}  \Lambda^\Delta_T.
\end{equation}  

\smallskip
Next, we use Propositions \ref{thm1} and \ref{prop:density_est}(a) to find an $\teps>0$ such that \eqref{loc_rate_conv} applies to all $T\le t\le T+\teps$ with this $\teps$. To conclude we argue that $\lim_{\Delta\downarrow0} \Lambda^\Delta_s=\Lambda_s$, $s\in(T,T+\teps]$. Indeed, for every $\eta>0$, there exists a $\overline{\Delta}_0>0$ with the property 
\begin{equation}\label{Delta1}
\sup_{s\in[t,t+\teps]} |\Lambda^{t;\Delta_0}_s-\Lambda_s|\le\frac{\eta}{2},\quad T\le t\le T+\teps,\quad \Delta_0\in(0,\overline{\Delta}_0]. 
\end{equation}
We now fix an $s\in(T,T+\teps]$. Then, $\lim_{\Delta\downarrow0} \lceil T/\Delta\rceil\Delta=T$, the right-continuity of $\Lambda$, and the preceding paragraph yield, for a sufficiently small $\overline{\Delta}\in(0,\overline{\Delta}_0)$,  
\begin{equation}\label{Delta2}
\lceil T/\Delta\rceil\Delta\le s,\,
\Lambda_{\lceil T/\Delta\rceil\Delta}-\Lambda_T\le  \frac{\eta}{4C_\alpha(\overline{\Delta}_0\!-\!\overline{\Delta},\teps)}
\,\text{and}\,\La_T - \La^\Delta_T \le \frac{\eta}{4C_\alpha(\overline{\Delta}_0\!-\!\overline{\Delta},\teps)} ,\;\Delta\!\in\!(0,\overline{\Delta}],
\end{equation}
where $C_\alpha(\overline{\Delta}_0-\overline{\Delta},\teps)\ge0$ is taken according to Corollary \ref{cor::discrete_compare} and is, in particular, non-increasing in the first argument.

\medskip

For any $\Delta\in(0,\overline{\Delta}]$ and integer multiple $\Delta_0\in[\overline{\Delta}_0-\overline{\Delta},\overline{\Delta}_0]$ of $\Delta$, we apply \eqref{Delta1}, followed by Corollary \ref{cor::discrete_compare} with $\CTT:=(0,\Delta,2\Delta,\ldots,\lceil T/\Delta\rceil\Delta,\lceil T/\Delta\rceil\Delta+\Delta_0,\lceil T/\Delta\rceil\Delta+2\Delta_0,\ldots)$
and \eqref{Delta2} to obtain
\begin{equation*}
\begin{split}
\La_s - \La^{\CTT}_s 
&= (\La_s - \La^{\lceil T/\Delta\rceil\Delta;\Delta_0}_s) 
+ (\La^{\lceil T/\Delta\rceil\Delta;\Delta_0}_s - \La^{\CTT}_s)  \\
&\le \frac{\eta}{2}+C_\alpha(\Delta_0,\teps)\cdot(\Lambda_{\lceil T/\Delta\rceil\Delta}-\La^\Delta_{\lceil T/\Delta\rceil\Delta}) \\
&\le\frac{\eta}{2}+C_\alpha(\overline{\Delta}_0-\overline{\Delta},\teps)\cdot\big((\Lambda_{\lceil T/\Delta\rceil\Delta}-\Lambda_T)+(\Lambda_T-\Lambda^\Delta_T)\big)\le\eta.
\end{split}
\end{equation*}
Thus, $\Lambda_s-\Lambda^\Delta_s\le \eta$ by Remark \ref{rmk::discrete_general_prelim} and, since $\eta>0$ was arbitrary, $\lim_{\Delta\downarrow0} \Lambda^\Delta_s=\Lambda_s$. \ep

\subsection{Extension beyond discontinuity points} \label{sec::special_case}

The main result of this subsection is the following analogue of Proposition \ref{pro::no_jump} for the case that $T\ge0$ is a discontinuity point of $\Lambda$.

\begin{prop}\label{prop_jump}
Suppose Assumption \ref{ass}. Let $T\ge0$ be a discontinuity point of $\Lambda$ such that $\Lambda^\Delta|_{[0,t]}\underset{\Delta\downarrow0}{\overset{\emph{M1}}{\longrightarrow}}\Lambda|_{[0,t]}$ for any continuity point $t\in(0,T)$ of $\Lambda$. Then,  $\Lambda^\Delta|_{[0,T+\teps]}\underset{\Delta\downarrow0}{\overset{\emph{M1}}{\longrightarrow}}\Lambda|_{[0,T+\teps]}$ for all $\teps>0$ small enough.
\end{prop}

The proof of Proposition \ref{prop_jump} builds on two lemmas. The first investigates the cumulative distribution function of $X^{\Delta}_{T+\theta}$ on a right neighborhood of $0$, for small $\Delta$, $\theta$ and on the event $\inf_{0\le s<T+\theta} X^\Delta_s>0$, and provides an estimate akin to the one in Proposition \ref{thm1}(ii). 

\begin{lem}\label{lem:jump-dens}
In the situation of Proposition \ref{prop_jump}, for any $\kappa>0$, there exist $\overline{\theta},\overline{\Delta}>0$ such that for all $\theta\in(0,\overline{\theta}]$ and all  $\Delta=(T+\theta)/1,\,(T+\theta)/2,\,\ldots\in(0,\overline{\Delta}]$,
\begin{equation}\label{eq::distrib_of_jump2}
\pp\Big(X^{\Delta}_{T + \theta} \le x,\,\inf_{0\le s<T+\theta} X_s^{\Delta} > 0\Big) \ge \frac{x-\kappa}{\alpha},\quad \kappa\le x\le\La_T - \La^{\Delta}_{T+\theta}. 
\end{equation}
\end{lem}

\noindent\textbf{Proof.}
By the minimality property \eqref{physical},
\[
\pp\Big( 0 < X_{T-} \le x,\,\inf_{0\le s<T} X_s > 0\Big) \ge \frac{x}{\alpha},\quad 	x \in [0, \La_T - \La_{T-}).
\]
For any $\Delta>0$, we have $X_{T-}=X^\Delta_T+\Lambda^\Delta_T-\Lambda_{T-}$ (since $X_{T-}+\Lambda_{T-}=X_{0-}+B_T=X^\Delta_T+\Lambda^\Delta_T$) and $X^\Delta\ge X$, so  
\[
\pp\Big( 0 < X^\Delta_T+\Lambda^\Delta_T-\Lambda_{T-} \le x,\,\inf_{0\le s < T} X^{\Delta}_s > 0\Big) \ge \frac{x}{\alpha},
\quad x \in [0, \La_T - \La_{T-}).
\]
We use this assertion upon setting $\ell=\Lambda^\Delta_T-\Lambda_{T-}$ to obtain for all $y\in[-\ell,\La_T - \La_{T-}-\ell)$,
\begin{eqnarray}
&& \pp\Big(X^\Delta_T \le y,\,\inf_{0\le s <T} X^{\Delta}_s > 0 \Big) 
= \pp\Big( X^\Delta_T+\ell \le y+\ell,\,\inf_{0\le s < T} X^{\Delta}_s > 0 \Big) 
\ge \frac{y+\ell}{\alpha},\quad\mathrm{and\;thus} \nonumber \\
&&\pp\Big(X^\Delta_T \le y,\,\inf_{0\le s <T} X^{\Delta}_s > 0 \Big)\ge\frac{\min(y+\ell,\Lambda_T-\Lambda_{T-})-\min(y+\ell,0)}{\alpha}, \quad y\in\rr. 
\label{eq::distrib_of_jump} 
\end{eqnarray}

\smallskip

Next, we turn to $X^\Delta_{T+\theta}$, let $\widetilde{x}:=x+\La^{\Delta}_{T+\theta}-\Lambda_{T-}$, and observe 
\begin{equation*}
\begin{split}
&\pp\Big(X^{\Delta}_{T + \theta} \le x,\,\inf_{0\le s < T} X_s^{\Delta} > 0\Big) 
= \pp\Big( X^{\Delta}_T + B_{T+\theta} - B_T - \La^{\Delta}_{T+\theta} + \La^{\Delta}_{T} \le x,\,\inf_{0\le s < T} X_s^{\Delta} > 0 \Big) \\
&\qquad\qquad\qquad\quad\;\;\;
= \frac{1}{\sqrt{2\pi}} \int_{-\infty}^{\infty} e^{-\frac{z^2}{2}}\, 
\pp\Big(X^{\Delta}_T - \La^{\Delta}_{T+\theta} +  \La^{\Delta}_T  \le x + \sqrt{\theta}z,\,\inf_{0\le s < T} X_s^{\Delta} > 0 \Big)\,\mathrm{d}z \\
&\qquad\qquad\qquad\quad\;\;\;
\ge \frac{1}{\sqrt{2\pi}\alpha} \int_{-\infty}^{\infty} e^{-\frac{z^2}{2}}\,\big(\min(\widetilde{x}+\sqrt{\theta}z,\Lambda_T-\Lambda_{T-})-\min(\widetilde{x}+\sqrt{\theta}z,0)\big)\,\mathrm{d}z \\
&\qquad\qquad\qquad\quad\;\;\;
\ge \frac{1}{\sqrt{2\pi}\alpha} \int_{-\widetilde{x}/\sqrt{\theta}}^{(\Lambda_T-\Lambda_{T-}-\widetilde{x})/\sqrt{\theta}} e^{-\frac{z^2}{2}}\,
(\widetilde{x}+\sqrt{\theta}z)\,\mathrm{d}z,
\end{split}
\end{equation*}
where we have used  \eqref{eq::distrib_of_jump} to get the first inequality. 
With Definition \ref{defn::discrete} and $(T+\theta)/\Delta\in\mathbb{N}$, 
\begin{equation*}
\pp\Big(X^{\Delta}_{T+\theta}\le x,\inf_{0\le s<T+\theta} X_s^{\Delta} > 0\Big) 
\ge\frac{1}{\sqrt{2\pi}\alpha} \int_{-\widetilde{x}/\sqrt{\theta}}^{(\Lambda_T-\Lambda_{T-}-\widetilde{x})/\sqrt{\theta}} 
\!\! e^{-\frac{z^2}{2}}(\widetilde{x}+\sqrt{\theta}z)\,\mathrm{d}z-\frac{\Lambda^\Delta_{T+\theta}\!
-\!\La^{\Delta}_{T-}}{\alpha}. 
\end{equation*}
Thus, given a $\kappa>0$, whenever $\theta>0$ is sufficiently small we have
\begin{equation}\label{numer_blowup}
\begin{split}
\pp\Big(X^\Delta_{T+\theta} \le x,\,\inf_{0\le s< T+\theta} X_s^{\Delta} > 0\Big) 
\ge \frac{x-\kappa/2}{\alpha}-\frac{\La_{T-} - \La^{\Delta}_{T-}}{\alpha},\\
\La_{T-} - \La^{\Delta}_{T+\theta}+\kappa/2\le x\le\La_T - \La^{\Delta}_{T+\theta}-\kappa,
\end{split}
\end{equation}
with the bounds on $x$ ensuring that the latter integration interval contains $[-\kappa/(2\sqrt{\theta}),\kappa/\sqrt{\theta}]$.  

\medskip		

We now claim that $\Lambda_{T-}-\Lambda^\Delta_{T-}\le\kappa/2$ for all $\Delta>0$ small enough. Indeed, this is obvious if $T=0$. For $T>0$, we let $t\in(0,T)$ be a continuity point of $\Lambda$ with the property $\Lambda_{T-}-\Lambda_t\le\kappa/4$. For all $\Delta>0$ small enough, $\Lambda_t-\Lambda^\Delta_t\le\kappa/4$, hence
\begin{equation*}
\Lambda_{T-}-\Lambda^\Delta_{T-}=(\Lambda_{T-}-\Lambda_t)+(\Lambda_t-\Lambda^\Delta_t)+(\Lambda^\Delta_t-\Lambda^\Delta_{T-})\le\frac{\kappa}{4}+\frac{\kappa}{4}+0=\frac{\kappa}{2}. 
\end{equation*}
In particular, $\La_{T-} - \La^{\Delta}_{T+\theta}\le\kappa/2$, and we deduce \eqref{eq::distrib_of_jump2} for $\kappa\le x\le\La_T - \La^{\Delta}_{T+\theta}-\kappa$ from \eqref{numer_blowup}.
For $\La_T - \La^{\Delta}_{T+\theta}-\kappa<x\le\La_T - \La^{\Delta}_{T+\theta}$, we use $\pp(X^{\Delta}_{T + \theta} \le x,\,\ldots\,) \ge \pp(X^{\Delta}_{T + \theta} \le x - \kappa,\,\ldots\,)$, apply \eqref{eq::distrib_of_jump2} with $x$ replaced by $x-\kappa$, and relabel $2\kappa$ as $\kappa$. \ep

\medskip

Next, we build on Lemma \ref{lem:jump-dens} to verify that the time-stepping scheme ``catches up'' with the physical solution. 

\begin{lem}\label{lem::special_case}
In the situation of Proposition \ref{prop_jump}, for any $\eta,\theta>0$, there exists some $\overline{\Delta}>0$ such that 
$\La_T - \La^\Delta_{T+\theta}\le\eta$, $\Delta\in(0,\overline{\Delta}]$.
\end{lem}
	
\noindent\textbf{Proof.} We establish the lemma by proving the following statement: If, for some $\eta> 0$, it holds that for all $\theta>0$ we can find some $\overline{\Delta}>0$ such that 
$\La_T - \La^{\Delta}_{T+\theta} \le \eta$, $\Delta\in(0,\overline{\Delta}]$, then the same holds for $\frac23\eta$ in place of $\eta$.	
To this end, we fix such an $\eta>0$ and any $\theta>0$. We seek a value of $\overline{\Delta}>0$ for which 
$\La_T - \La^\Delta_{T+\theta} \le \frac23 \eta$, $\Delta\in(0,\overline{\Delta}]$. For any $\theta_0>0$, 
the hypothesis of the statement yields a $\overline{\Delta}_0>0$ such that 
$\La_T - \La^{\Delta_0}_{T+\theta_0} \le \eta$, $\Delta_0\in(0,\overline{\Delta}_0]$. 
By Lemma \ref{lem:jump-dens}, we may select $\theta_0$ and $\Delta_0$ to which \eqref{eq::distrib_of_jump2} also applies, with a $\kappa=\kappa(\eta,\theta)>0$ to be determined below. Our aim is to 
show that if $\theta_0$ is picked from a suitable interval $(0,\overline{\theta}_0]$, then $\widetilde{\Delta}>0$ can be chosen such that for all $\Delta=\Delta_0/1,\,\Delta_0/2,\,\ldots\in(0,\widetilde{\Delta}]$,  
\begin{eqnarray}\label{taueqn}
\La_T - \La^{\CTT}_{T+\theta} \le \frac23 \eta,\;\text{with}\;
\CTT := (0,\Delta_0,2\Delta_0,\ldots, T+\theta_0,T+\theta_0+\Delta,T+\theta_0+2\Delta,\dots).
\end{eqnarray}
The claim then follows, since $\Lambda^\Delta_{T+\theta}\ge\La^{\CTT}_{T+\theta}$ by Remark \ref{rmk::discrete_general_prelim} (thus, $\Lambda_T-\Lambda^\Delta_{T+\theta}\le \frac23 \eta$) and the range $\bigcup_{\theta_0\in(0,\overline{\theta}_0]} \{(T+\theta_0)/1,(T+\theta_0)/2,\ldots\}\cap(0,\overline{\Delta}_0]$ of suitable $\Delta_0$ contains the interval $(0,\min(\overline{\theta}_0,\overline{\Delta}_0)]$ allowing us to take  $\overline{\Delta}=\min(\widetilde{\Delta},\overline{\theta}_0,\overline{\Delta}_0)>0$.  

\medskip
		
Fix any $\Delta=\Delta_0/1,\,\Delta_0/2,\,\ldots$ and
let $\ell_n:=\La^{\CTT}_{T+\theta_0+n\Delta}-\La^{\Delta_0}_{T+\theta_0}$, $n=0,\,1,\,\ldots\,$. 
Then,
\[
\begin{split}
\ell_{n+1} &= \alpha\pp\Big(\inf_{0\le s<T+\theta_0} X_s^{\Delta_0} > 0,\,\min_{0\le m \le n} X^{\CTT}_{T +\theta_0+m\Delta} < 0\Big) \\
& \ge 
\alpha\pp\Big(\inf_{0\le s<T+\theta_0} X_s^{\Delta_0} > 0,\,X^{\CTT}_{T+\theta_0+n\Delta} < 0\Big) \\
& = \frac{\alpha}{\sqrt{2\pi}} \int_{-\infty}^{\infty} e^{-\frac{z^2}{2}}\,
\pp\Big(\inf_{0\le s<T+\theta_0} X_s^{\Delta_0} > 0,\,X^{\Delta_0}_{T+\theta_0} < \sqrt{n\Delta} z + \ell_n \Big)\,\mathrm{d}z.
\end{split}
\]
Next, we set $\widetilde{\ell}=\Lambda_T-\Lambda^{\Delta_0}_{T+\theta_0}$ and apply 
\eqref{eq::distrib_of_jump2} in the form
\[
\pp\Big(X^{\Delta_0}_{T+\theta_0} < x,\,\inf_{0\le s<T+\theta_0} X_s^{\Delta_0} > 0\Big) 
\ge\frac{\min(x-\kappa,\widetilde{\ell}-\kappa)-\min(x-\kappa,0)}{\alpha},\quad x\in\rr
\]
to obtain the explicit recursive inequality 
\[
\ell_{n+1}\ge\frac{1}{\sqrt{2\pi}} \int_{-\infty}^{\infty} e^{-\frac{z^2}{2}}\,
\Big(\!\min\big(\sqrt{n\Delta}z + \ell_n - \kappa,\widetilde{\ell}-\kappa\big) 
- \min\big(\sqrt{n\Delta}z + \ell_n - \kappa,0\big)\!\Big)\,\mathrm{d}z.
\]
We rewrite the latter using the standard Gaussian cumulative distribution function $\Phi$:
\[
\ell_{n+1}\ge\frac{1}{\sqrt{2\pi}} \int_{(\kappa-\ell_n)/\sqrt{n\Delta}}^{(\widetilde{\ell} - \ell_n)/\sqrt{n\Delta}}  e^{-\frac{z^2}{2}}\,\big(\sqrt{n\Delta}z+\ell_n-\kappa\big)\,\mathrm{d}z 
		 + (\widetilde{\ell}-\kappa)\bigg(1-\Phi \bigg(\frac{\widetilde{\ell}-\ell_n}{\sqrt{n\Delta}}\bigg)\!\bigg).
\]
Integrating by parts and introducing 
\[
f(x):=\int_{-\infty}^x \Phi(z)\,\mathrm{d}z=x\Phi(x) + \frac{1}{\sqrt{2\pi}}\,e^{-\frac{x^2}{2}},\quad x\in\rr
\]
we further deduce 
\[
\ell_{n+1} \ge \widetilde{\ell} - \kappa - \sqrt{n\Delta}\,\bigg(f\bigg(\frac{\widetilde{\ell}-\ell_n}{\sqrt{n\Delta}} \bigg) - f\bigg(\frac{\kappa-\ell_n}{\sqrt{n\Delta}} \bigg)\!\bigg).
\]
In view of $\Phi(x)+\Phi(-x)=1$, $x\in\rr$, we have $f(x)=f(-x)+x$, $x\in\rr$, so 
\[
\ell_{n+1} \ge \ell_n - \kappa + \sqrt{n\Delta}\,\bigg(f\bigg(\frac{\kappa-\ell_n}{\sqrt{n\Delta}}\bigg) 
-f\bigg(\!-\frac{\widetilde{\ell}-\ell_n}{\sqrt{n\Delta}}\bigg)\!\bigg).
\]
Since $f'=\Phi$ is increasing, we arrive at
\begin{equation}\label{eq::jump_recursion}
\ell_{n+1} \ge \ell_n - \kappa + \Phi\bigg(\!-\frac{\widetilde{\ell}-\ell_n}{\sqrt{n\Delta}} \bigg)
\cdot(\kappa-2\ell_n+\widetilde{\ell}),
\end{equation}
provided that $\kappa-2\ell_n+\widetilde{\ell}\ge 0$.
		

\medskip

To prove \eqref{taueqn}, and thus the lemma, we assume $\widetilde{\ell}>2\eta/3$ (otherwise \eqref{taueqn} holds with $\theta_0=\theta$ and $\widetilde{\Delta}=\overline{\Delta}_0$) and show that if $\theta_0\in(0,\theta)$ is taken small enough and the integers $(T+\theta_0)/\Delta_0$ and $\Delta_0/\Delta$ are made large enough, then one 
can find an $n\ge0$ with $\theta_0+n\Delta\le\theta$ and $\widetilde{\ell}/3\le \ell_n$. The inequality $\La_T - \La^{\CTT}_{T+\theta} \le  2\eta/3$ in \eqref{taueqn} then follows from 
\[
\La_T-\La^{\CTT}_{T+\theta_0+n\Delta}
=(\La_T-\La^{\Delta_0}_{T+\theta_0})+(\La^{\Delta_0}_{T+\theta_0}-\La^{\CTT}_{T+\theta_0+n\Delta})
\le \frac{2(\La_T-\La^{\Delta_0}_{T+\theta_0})}{3}\le\frac{2\eta}{3}.
\]
Arguing by contradiction, suppose that $\ell_n<\widetilde{\ell}/3$ for all $n=0,\,1,\,\ldots,\,\lfloor(\theta - \theta_0)/\Delta\rfloor=:n_0$. Summing \eqref{eq::jump_recursion} over $0\le n\le n_0-1$ upon noting that $\kappa-2\ell_n+\widetilde{\ell}>\kappa+\widetilde{\ell}/3>\kappa+\ell_0=\kappa>0$ for such $n$ we get
\begin{eqnarray*}
\ell_{n_0} \ge -\kappa n_0 + \sum_{n = 0}^{n_0-1} \Phi\bigg(\!- \frac{\widetilde{\ell}-\ell_n}{\sqrt{n\Delta}}\bigg)
\cdot(\kappa-2\ell_n+\widetilde{\ell}).
\end{eqnarray*}
Using $\kappa-2\ell_n+\widetilde{\ell}>\widetilde{\ell}/3$ and $\widetilde{\ell}-\ell_n\le\widetilde{\ell}$ we deduce 
\[
\begin{split}
\ell_{n_0} 
\ge -\kappa n_0 
+\frac{\widetilde{\ell}}{3}\,\sum_{n = 0}^{n_0-1} \Phi\bigg(\!-\frac{\widetilde{\ell}}{\sqrt{n\Delta}} \bigg) 
&\ge n_0\bigg(\!-\kappa+\frac{\widetilde{\ell}}{6}\, \Phi\bigg(\!-\frac{\sqrt2\widetilde{\ell}}{\sqrt{\theta-\theta_0-2\Delta}} \bigg)\!\bigg) \\
&\ge n_0\bigg(\!-\kappa+\frac{\eta}{9}\, \Phi\bigg(\!-\frac{\sqrt2\eta}{\sqrt{\theta-\theta_0-2\Delta}} \bigg)\!\bigg),
\end{split}
\]
where we have dropped the summands with $n<n_0-\lceil n_0/2\rceil$ and have employed a common lower bound for the remaining summands. We choose $\kappa=\kappa(\eta,\theta)>0$ such that the latter bracket is positive, observing that the smaller $\theta_0+2\Delta>0$ is, the larger this bracket becomes. Lastly, we pick $\Delta>0$ dividing $\Delta_0$ so that $n_0$ is large and the final lower bound on $\ell_{n_0}$ is at least $\Lambda_T/3\ge\widetilde{\ell}/3$, giving us the desired contradiction. The proof is complete. \ep
		
\medskip	
We are now ready to show Proposition \ref{prop_jump}.

\medskip

\noindent\textbf{Proof of Proposition \ref{prop_jump}.} We start by using Propositions \ref{thm1} and \ref{prop:density_est}(a) to find an $\teps>0$ such that \eqref{loc_rate_conv} applies to all $T\le t\le T+\teps$ with this $\teps$. In particular, upon fixing $s\in(T,T+\teps]$ and $\eta>0$ we can select some $\overline{\Delta}_0>0$ so that
\begin{equation}\label{eta/2bnd}
\Lambda_s-\Lambda^{t;\Delta_0}_s\le\frac{\eta}{2},\quad t\in[T,s],\quad\Delta_0\in(0,\overline{\Delta}_0]. 
\end{equation}
Next, we recall the quantity $C_\alpha(\overline{\Delta}_0/2,s-T)\ge0$ from Corollary \ref{cor::discrete_compare} and choose $\theta_0\in(0,s-T]$ with the property
\begin{equation}\label{rightcont}
\Lambda_{T+\theta_0}-\Lambda_T\le\frac{\eta}{4C_\alpha(\overline{\Delta}_0/2,s-T)}. 
\end{equation} 
Finally, we rely on Lemma \ref{lem::special_case} to pick $\widetilde{\Delta}>0$ so that 
\begin{equation}\label{catchup}
\Lambda_T-\Lambda^\Delta_{T+\theta_0/2}\le\frac{\eta}{4C_\alpha(\overline{\Delta}_0/2,s-T)},\quad\Delta\in(0,\widetilde{\Delta}]. 
\end{equation}
Then, for $\Delta\in(0,\min(\widetilde{\Delta},\overline{\Delta}_0/2,\theta_0/2)]$ and $\Delta_0\in[\overline{\Delta}_0/2,\overline{\Delta}_0]$ such that $\Delta_0/\Delta$ is an integer, we get with $t=\lfloor (T+\theta_0)/\Delta\rfloor\Delta\in[T+\theta_0/2,T+\theta_0]$ and $\CTT=(0,\Delta,2\Delta,\ldots,t,t+\Delta_0,t+2\Delta_0,\ldots)$, 
\begin{equation}
\begin{split}
\Lambda_s-\Lambda^\Delta_s\le \Lambda_s-\Lambda^{\CTT}_s
&=(\Lambda_s-\Lambda^{t;\Delta_0}_s)+(\Lambda^{t;\Delta_0}_s-\Lambda^{\CTT}_s) \\
&\le\frac{\eta}{2}+C_\alpha(\overline{\Delta}_0/2,s-T)(\Lambda^{t;\Delta_0}_t-\Lambda^{\CTT}_t) \\
&=\frac{\eta}{2}+C_\alpha(\overline{\Delta}_0/2,s-T)\big((\Lambda_t-\Lambda_T)
+(\Lambda_T-\Lambda^\Delta_t)\big)\le\eta,
\end{split}
\end{equation}
where we have used Remark \ref{rmk::discrete_general_prelim}, \eqref{eta/2bnd}, Corollary \ref{cor::discrete_compare}, Definition \ref{defn::discrete_general}, \eqref{rightcont}, and \eqref{catchup} in this order. 
Since $s\in(T,T+\teps]$ and $\eta>0$ were arbitrary, we conclude that $\lim_{\Delta\downarrow0} \Lambda^\Delta_s=\Lambda_s$, $s\in(T,T+\teps]$, thus obtaining the proposition. \ep

\subsection{Proof of Theorem \ref{thm::main}} \label{sec::post_jump}

Theorem \ref{thm::main} follows directly from Propositions \ref{pro::no_jump}~and~\ref{prop_jump}. Indeed, we can argue by contradiction and assume the existence of continuity points $s\in(0,\infty)$ of $\Lambda$ for which the convergence $\Lambda^\Delta|_{[0,s]}\underset{\Delta\downarrow0}{\overset{\text{M1}}{\longrightarrow}}\Lambda|_{[0,s]}$ does not hold. Let $T\ge0$ be the infimum of such continuity points. If $T$ is a continuity point of $\Lambda$, we rely on Proposition~\ref{pro::no_jump} to infer $\Lambda^\Delta|_{[0,T+\teps]}\underset{\Delta\downarrow0}{\overset{\text{M1}}{\longrightarrow}}\Lambda|_{[0,T+\teps]}$ for all $\teps\ge0$ small enough. If $T$ is a discontinuity point of~$\Lambda$, we apply Proposition \ref{prop_jump} to get $\Lambda^\Delta|_{[0,T+\teps]}\underset{\Delta\downarrow0}{\overset{\text{M1}}{\longrightarrow}}\Lambda|_{[0,T+\teps]}$ for all $\teps>0$ small enough. This contradiction to the definition of $T$ yields the theorem. \ep 

\begin{rmk}
It is worth noting that the proof of Theorem \ref{thm::main} remains intact if instead of Assumption \ref{ass}(a), the boundedness of $f$ on $[0,\infty)$ together with the conclusions of Lemma~\ref{lemma:a priori} and Corollary \ref{cor:mod_of_cont} are assumed and $\psi:=\Psi':\,(0,\delta]\to\big(0,\frac{1}{\alpha}\big]$ is a strictly increasing function. Indeed, the proof of Theorem \ref{thm::main} can be then repeated word by word, and \eqref{loc_rate_conv}, used in the proofs of Propositions \ref{pro::no_jump} and \ref{prop_jump}, can be obtained via Remark \ref{alt_ass}. 
\end{rmk}




\section{Numerical simulations}  \label{sec::numerical}

In this last section, we examine the convergence of the time-stepping scheme (recall Definition \ref{defn::discrete}) numerically, for various initial densities $f$ and resulting functions $\Lambda$ with and without discontinuities. We write $L^\Delta$ for $\Lambda^\Delta/\alpha$, so that $L^\Delta_{n\Delta}
= \pp\big(\!\min_{0\le m\le n-1} X_{m\Delta}^\Delta<0\big)$, $n=1,\,2,\,\ldots\,$, 
and simulate the latter probabilities by a Monte Carlo particle method with $N=10^7$ particles, following \cite[Algorithm 1]{KR}. This fully implementable scheme
is given for completeness in the appendix, where we also derive the convergence as $N\rightarrow \infty$.
This extends the convergence result from \cite{KR}, which is restricted to the regular case studied there, and additionally shows the order $O(1/\sqrt{N})$ for fixed $\Delta$.

\subsection{Initial density vanishing at zero and no discontinuity}\label{sim1/2}

We first consider an example with $\lim_{x\downarrow0} f(x)=0$. 
To this end, we let $X_{0-}$ be $\Gamma(3/2,1/2)$-distributed and note that $f(x)\le Cx^{1/2}$, $x\ge0$ for a suitable constant $C<\infty$. Further we fix the time interval $[0,T]=[0,0.8]$ and $\alpha = 1.3$.

In Figure 1(a), we exhibit $L^{0.8/n}$ on the time interval $[0,0.8]$ for different numbers $n$ of time steps (or equivalently different values of $\Delta=0.8/n$). The numerical simulations indicate that there is no discontinuity in this setting.

\captionsetup[subfigure]{labelfont=rm}

\begin{figure}[h]
	\begin{center}
		\subfloat[]{\includegraphics[width=0.5\textwidth]{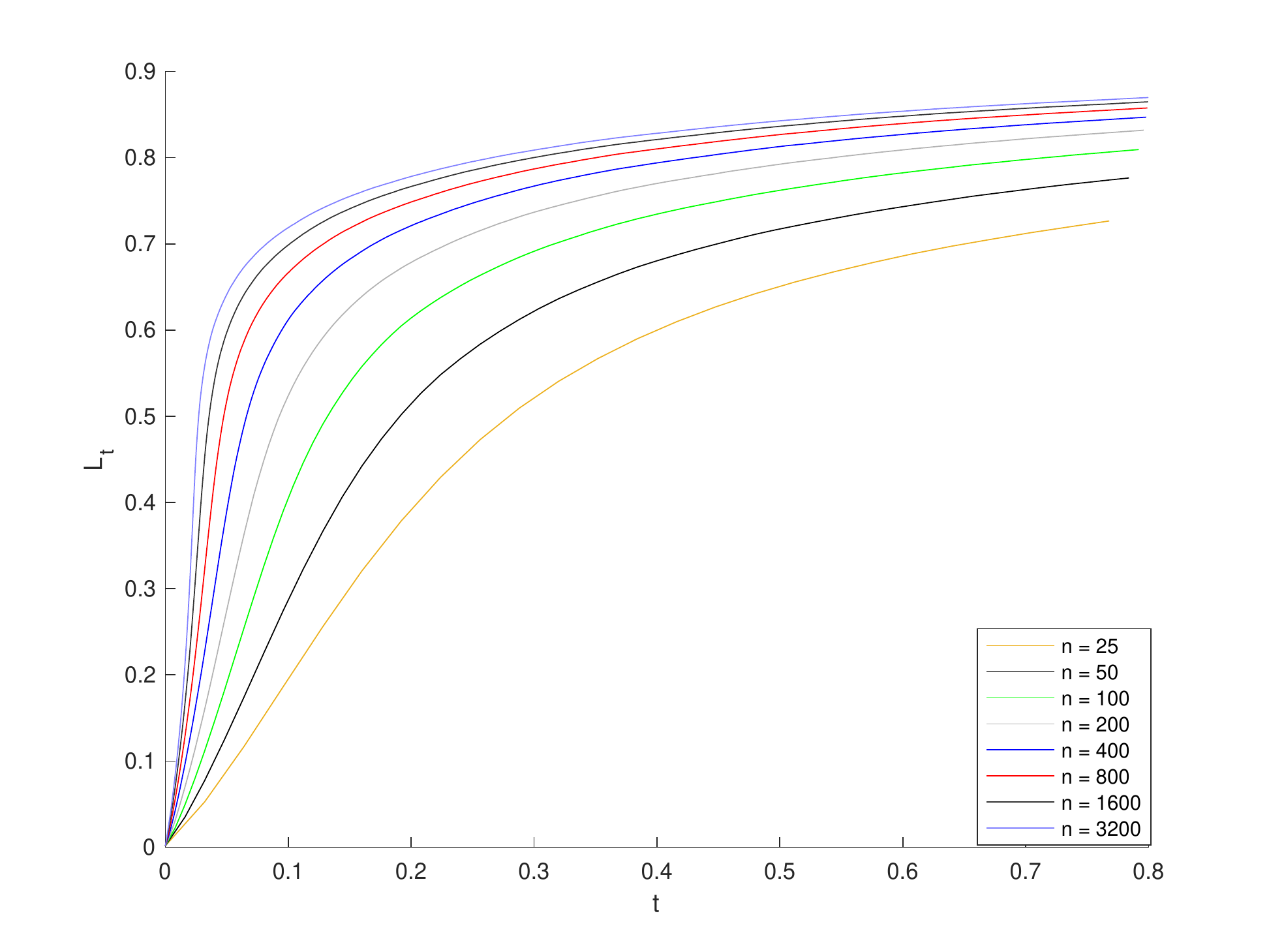}\label{fig:example1a}}
		\subfloat[]{\includegraphics[width=0.5\textwidth]{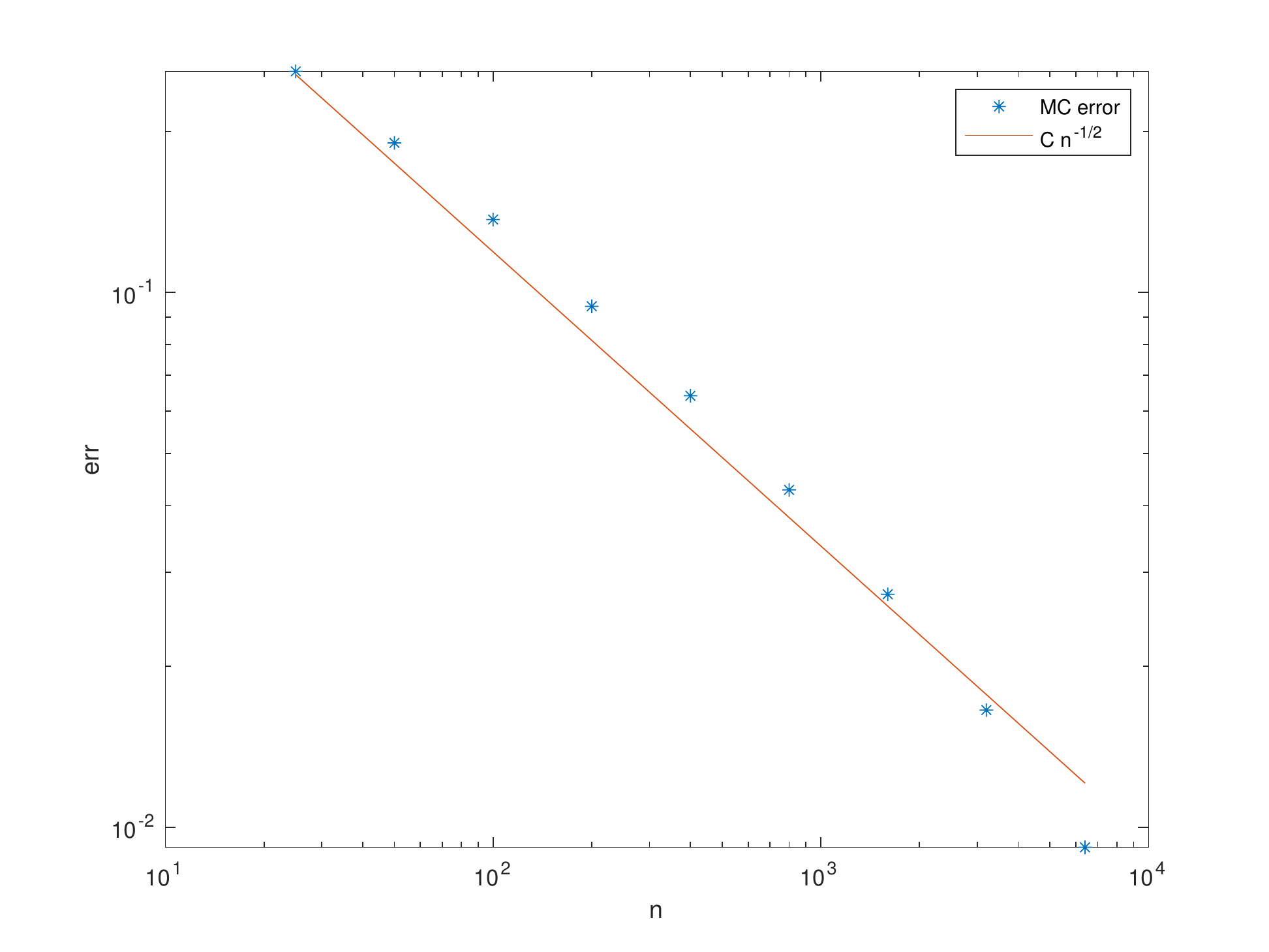}\label{fig:example1b}} \\
	\end{center}
	\caption{Example with $\lim_{x\downarrow0} f(x)=0$ and no discontinuity for different values of $n$: (a) absorption probabilities $L^{0.8/n}_t=\Lambda^{0.8/n}_t/1.3$, $t\in[0,0.8]$, and (b) error in the supremum norm over the time interval $[0,0.8]$.}
 	\label{fig:example1}
\end{figure}	

We approximate the error in the supremum norm by 
$\sup_{t\in [0,0.8]} |L_t^{0.8/n} - \bar{L}_t|$,
where a reference solution $\bar{L}:=L^{0.8/\bar{n}}$ is computed on a refined time mesh with $\bar{n} = 25600$ points.
The convergence rate suggested by Figure 1(b) is around 0.5, as confirmed visually by a comparison with a reference line (on the log-scale) of slope $-0.5$. This is expected from Remark \ref{rm_order}(a) and has been also observed in \cite[Section 4.2]{KR} for the pointwise error. We note that the two rightmost data points deviate slightly from the predicted line due to the comparison with a reference solution computed on a refined time mesh. Accordingly, the convergence order estimated by regression is slightly biased high, namely 0.55.

\subsection{Initial density of $1/\alpha$ at zero and no discontinuity}\label{sec5.3}

We turn to the numerical examination of the rate of local convergence when $\lim_{x\downarrow0} f(x)=\frac{1}{\alpha}$ (cf.~Theorem \ref{thm:loc_conv})
and illustrate that the rate may indeed be arbitrarily low. 
To this end, we pick the initial density
\begin{equation}
\label{dens_a_alph}
	f(x) = \left\{
	\begin{tabular}{ll}
		$\frac{1}{\alpha} - c x^{a}$, & \quad $0 \le x \le A$, \\
		$0$, & \quad $x > A$,
	\end{tabular} \right.
\end{equation}
for $\alpha>0$ and $a>0$, which we vary in the tests, and where $A>0$ is determined by $\int_0^\infty f(x)\,\mathrm{d}x=1$ for given $c>0$, the latter being sufficiently small. This choice corresponds to taking $\psi(x)=c x^a$ and $\Psi(x)=c x^{a+1}/(a+1)$ in Theorem \ref{thm:loc_conv}.
Moreover, we let $T=10^{-4}$ be small enough to precede a possible discontinuity. 

\medskip

\begin{figure}[h]
	\begin{center}
		\subfloat[]{\includegraphics[width=0.5\textwidth, height=0.3\textwidth]{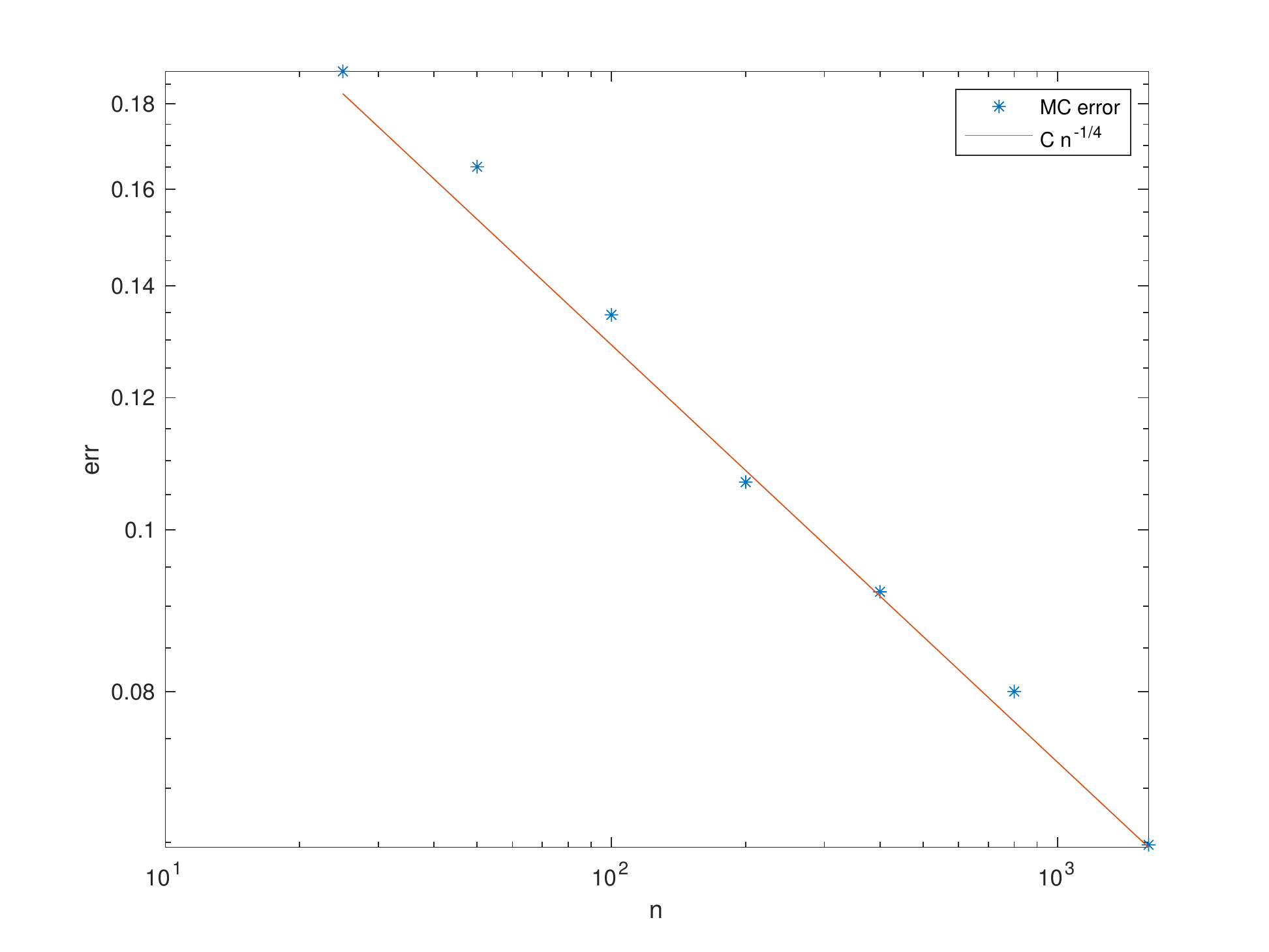}\label{fig:a1alpha05}}
		\subfloat[]{\includegraphics[width=0.5\textwidth, height=0.3\textwidth]{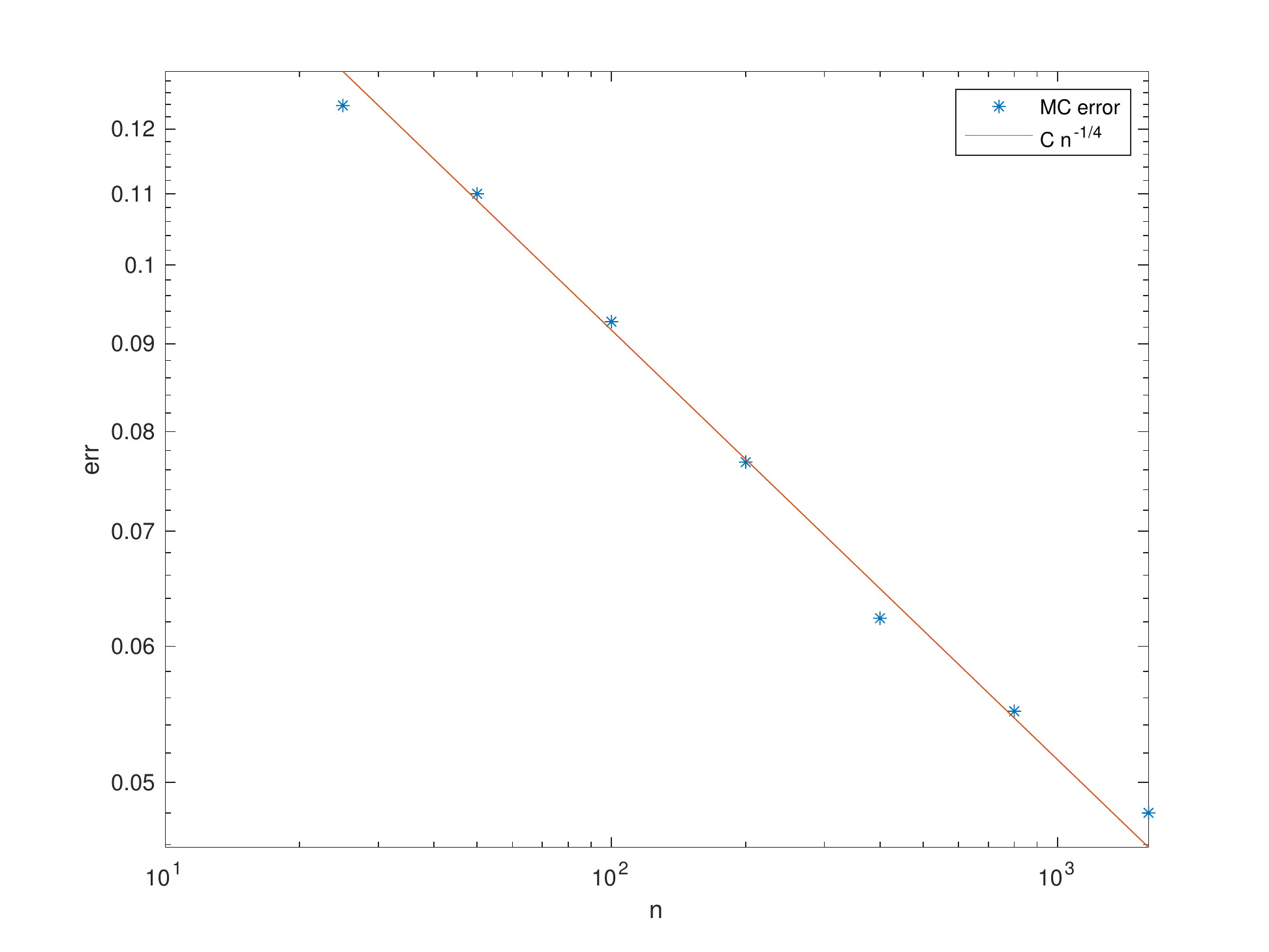}\label{fig:a1alpha1}} \\
		\subfloat[]{\includegraphics[width=0.5\textwidth, height=0.3\textwidth]{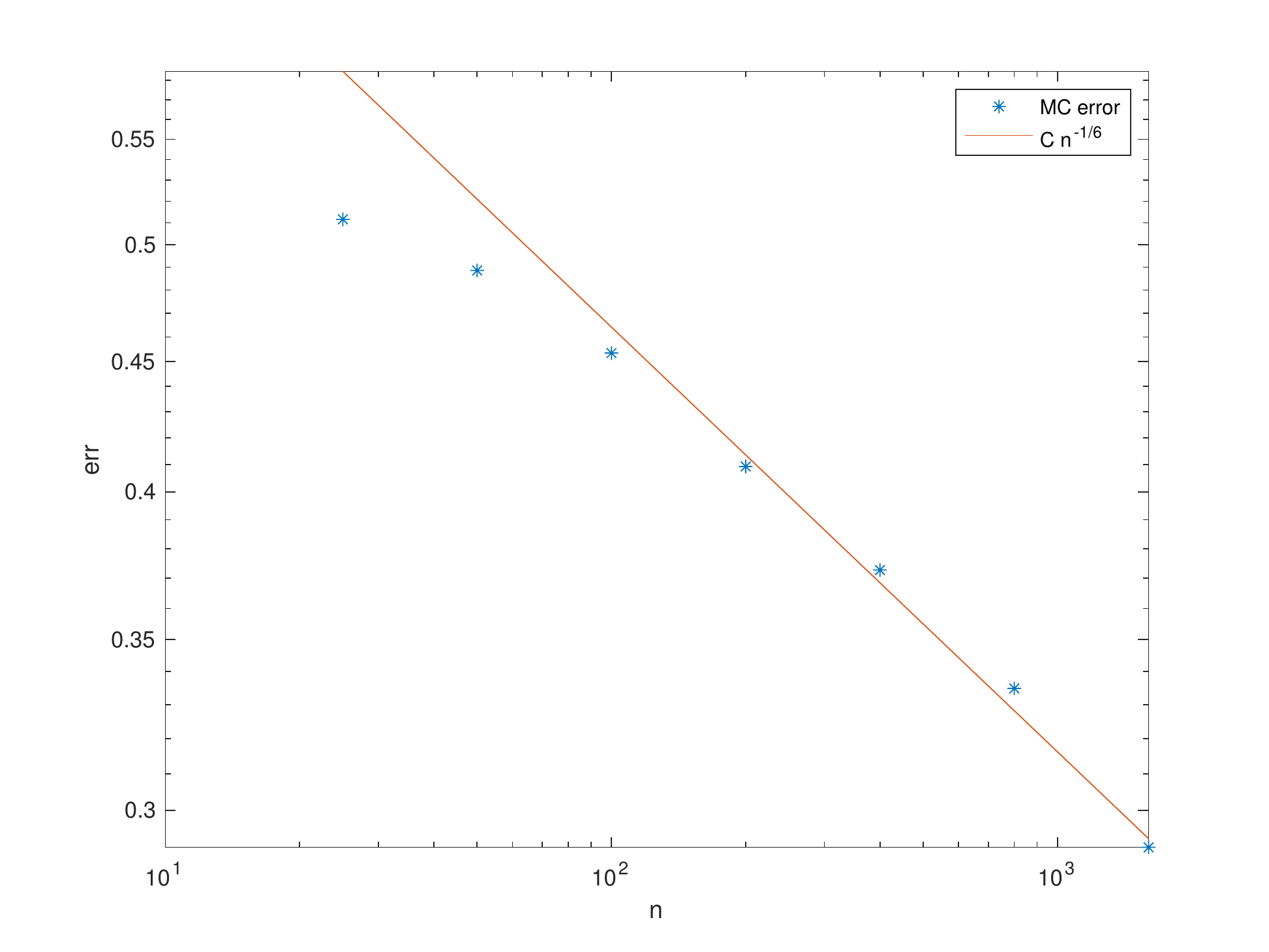}\label{fig:a2alpha05}}
		\subfloat[]{\includegraphics[width=0.5\textwidth, height=0.3\textwidth]{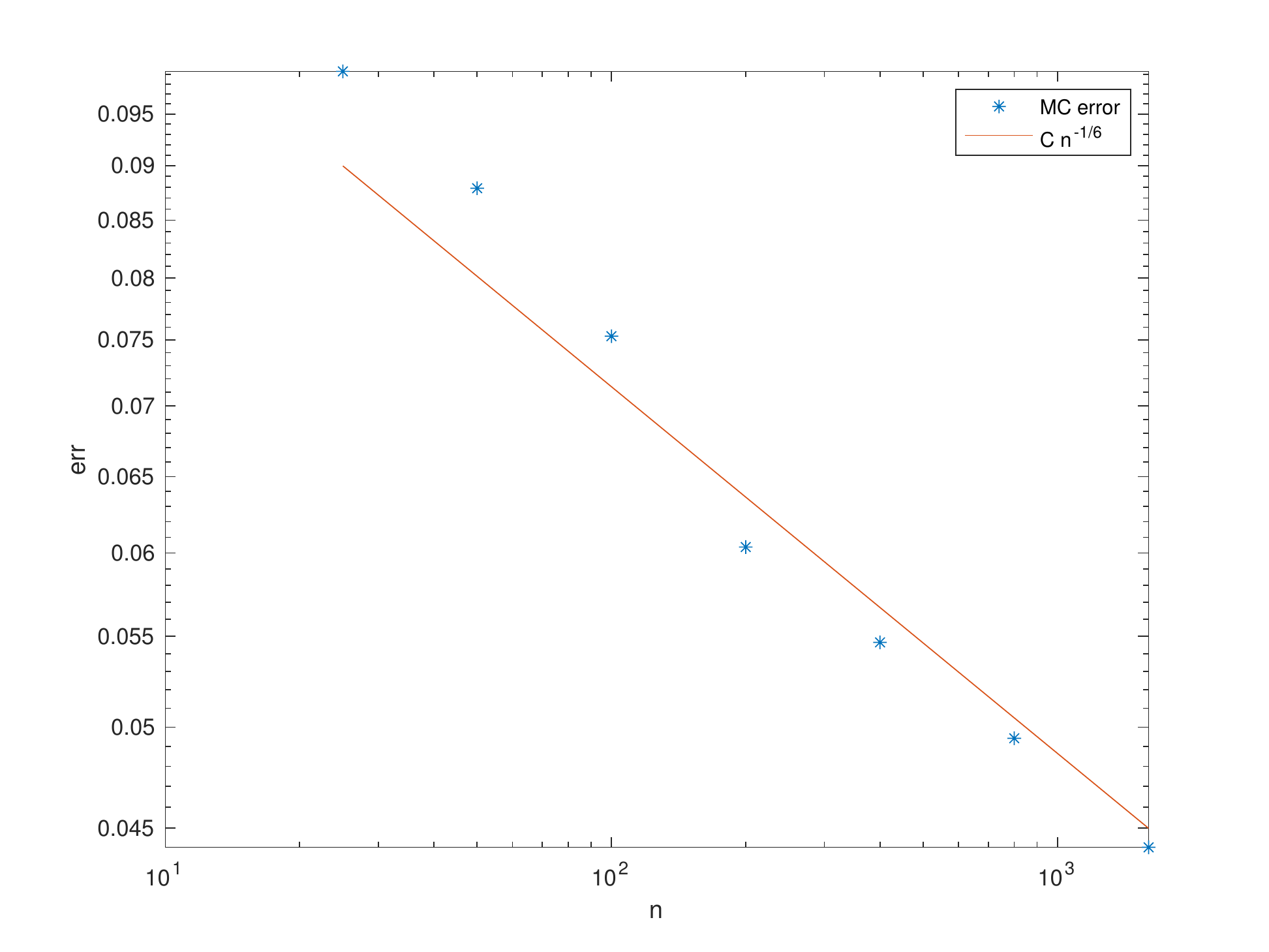}\label{fig:a2alpha1}} \\
		\subfloat[]{\includegraphics[width=0.5\textwidth, height=0.3\textwidth]{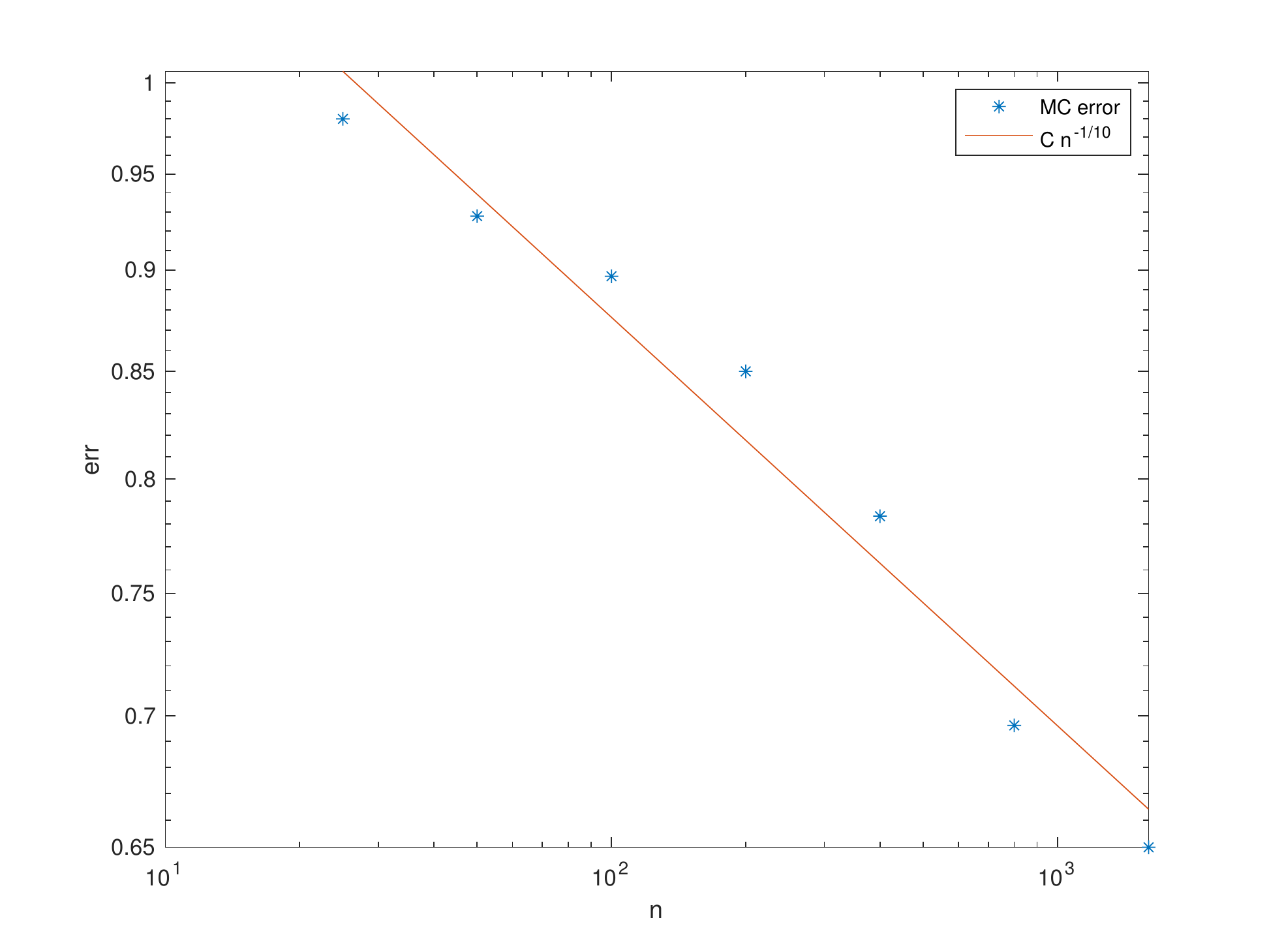}\label{fig:a4alpha05}}
		\subfloat[]{\includegraphics[width=0.5\textwidth, height=0.3\textwidth]{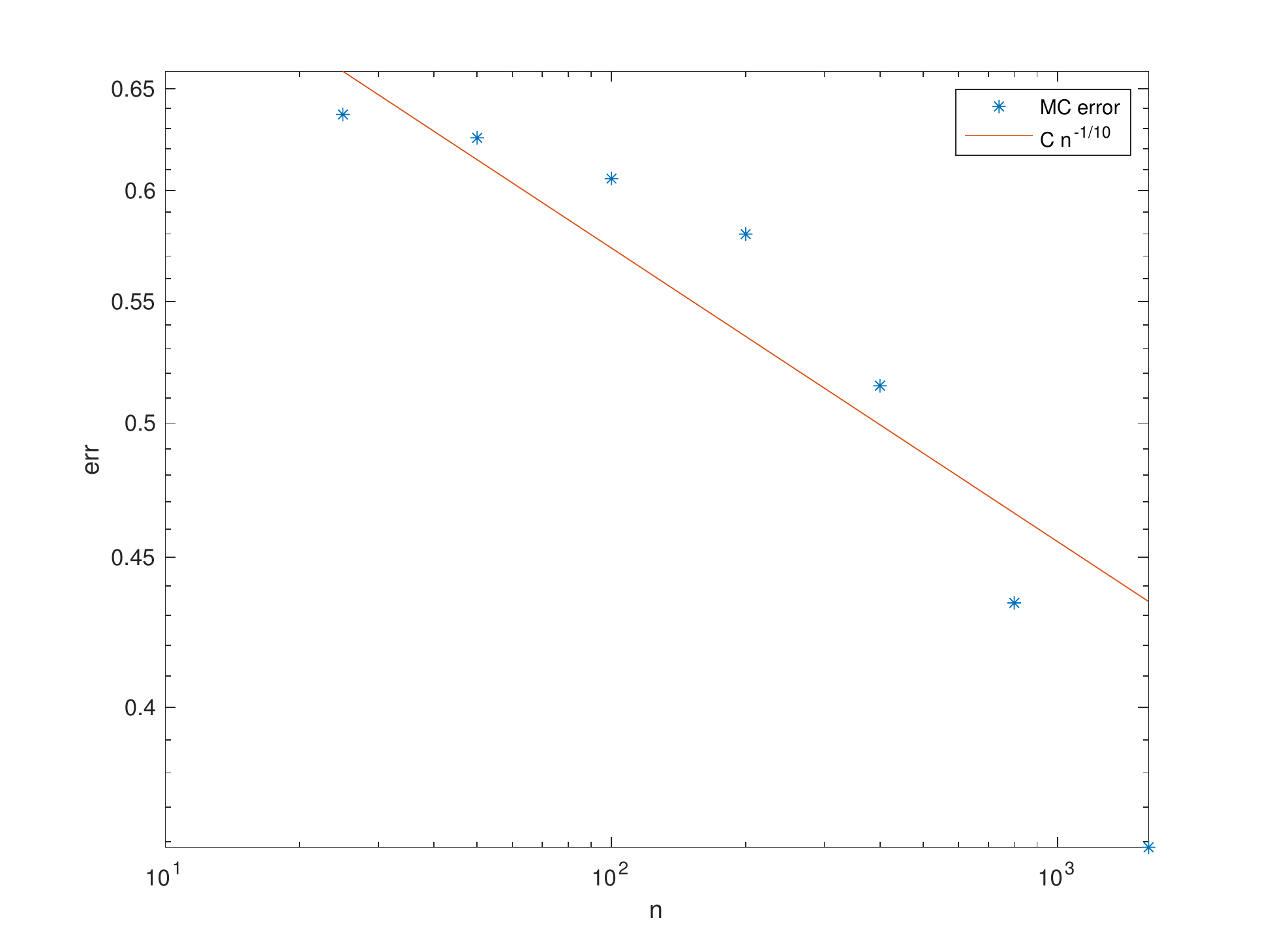}\label{fig:a4alpha1}}
	\end{center}
	\caption{Error in the supermum norm over $[0,10^{-4}]$ for the initial densities of \eqref{dens_a_alph} with: (a) $a=1$, $\alpha=0.5$, (b) $a=1$, $\alpha=1$, (c) $a=2$, $\alpha=0.5$, (d) $a=2$, $\alpha=1$, (e) $a=4$, $\alpha=0.5$, and (f) $a=4$, $\alpha=1$. We showcase the error for different $n$, and a reference line of slope $-1/(2(a+1))$.
	}
 	\label{fig:example2new}
\end{figure}

To exhibit the rate of convergence we plot in Figure \ref{fig:example2new} the error in the supremum norm as a function of the number $n$ of time steps on a log-log scale, for $\alpha\in\{0.5,1\}$ and $a\in\{1, 2, 4\}$. (We also adjust the constant $c$ for each $\alpha$ and $a$.)
These errors are compared to the reference line of slope $-1/(2(a+1))$, which reflects the term on the second line in \eqref{roc}.
The data for most parameter pairs suggest a convergence of polynomial order consistent with that term. As already pointed out in Remark \ref{rm_order}(b), the asymptotically dominant term on the right-hand side of \eqref{roc} is the first, logarithmic one, indicating that the theoretical result may not be sharp. The rate appears independent of $\alpha$, but the constant of proportionality of the error seems to be affected by it when one contrasts the errors in the respective left and right plots. The last two plots, for $a=4$ and $\alpha \in \{0.5,1 \}$, are less well explained by this. We find the convergence for large $a$ to be extremely slow, which may explain why the true asymptotic order is not detected for computationally feasible numbers of particles and time steps. 

\medskip

We also show in Table \ref{table:convergence_rate}, for varying values of $\alpha$ and $a$, the convergence rate estimated by regression from the logarithmic errors in the supremum norm for $n=25,\,50,\,\ldots,\,800$. The estimated rates are mostly independent of $\alpha$, and in many cases close to $1/(2(a+1))$, this being the rate suggested by the term on the second line in \eqref{roc}. Some significant deviations are again found for larger $a$ and $\alpha$, which may again be due to the large errors in such cases.
\begin{table}[H]
	\begin{center}
	\begin{tabular}{| c | c | c | c | c | c | c |}
  \hline			
$a \backslash \alpha$ & 1/(2($a$+1)) &   0.1 & 0.5 & 1.0 & 1.5\\
 \hline
1 & 1/4 =  0.25 & 0.234  & 0.259 & 0.261 & 0.261 \\
2 & 1/6 $\approx$ 0.17 & 0.169 &0.175 & 0.180 & 0.166  \\
4 & \!\!\!\!\!\!\!\! 1/10 = 0.1 &  0.104 & 0.125 & 0.101 & 0.089 \\
8 & 1/18 $\approx$ 0.056 & 0.053 & 0.049 & 0.025 & 0.033 \\ \hline
\end{tabular}
\vspace{15pt}
\caption{Estimated convergence rates on the time interval $[0,10^{-4}]$ for different values of $a$ and $\alpha$.}
\vspace{-20pt}
\label{table:convergence_rate}
\end{center}
\end{table}
%
%
%

\subsection{ 
Example with discontinuity} 

\begin{figure}[h]
	\begin{center}
		\includegraphics[width=0.7\textwidth,height=0.45\textwidth]{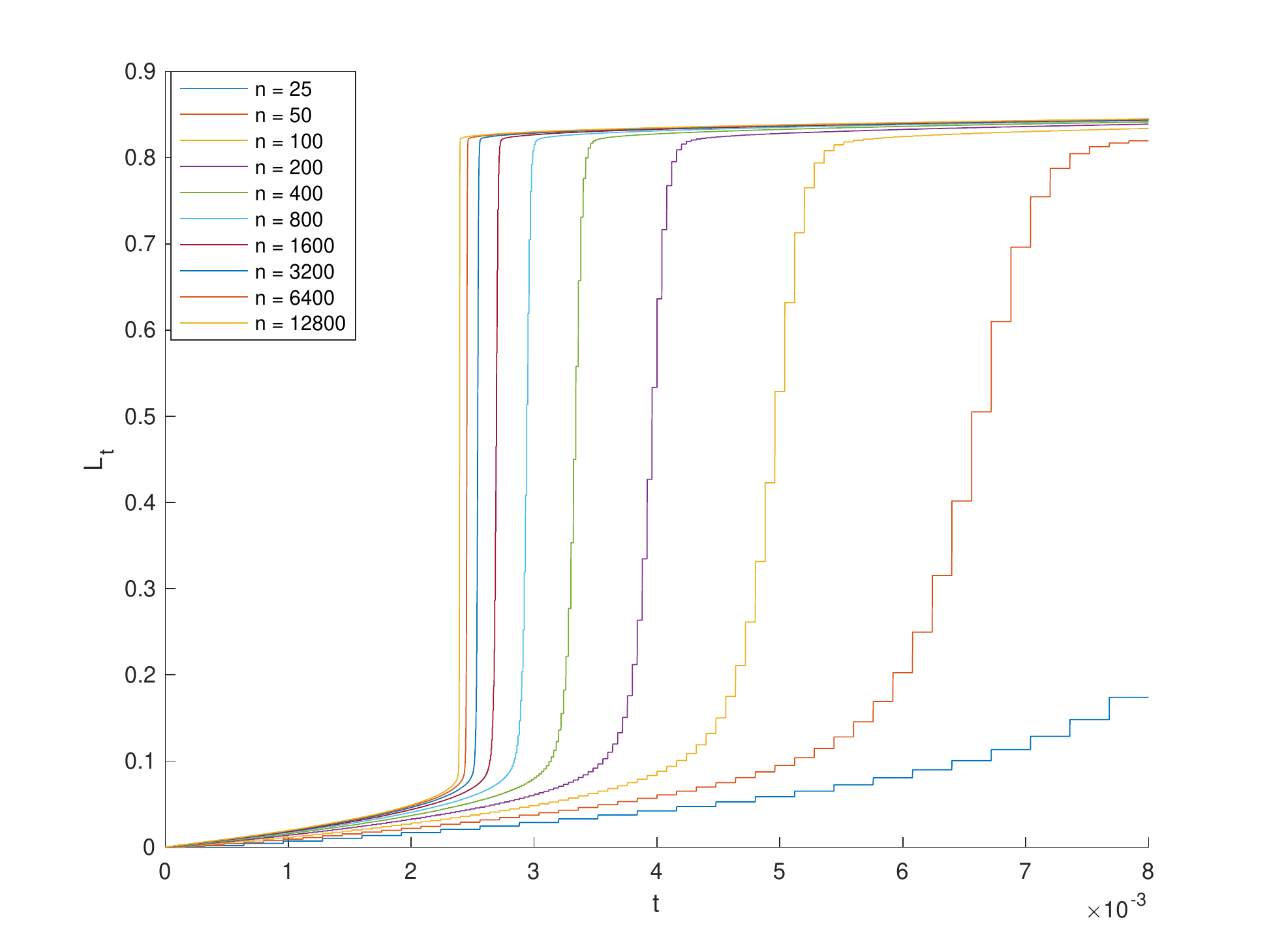}
	\end{center}
	\vspace{-17pt}
	\caption{Absorption probabilities $L^{0.008/n}_t=\Lambda^{0.008/n}_t/1.5$ for varying $n$.}
	\label{fig:example2}
\end{figure}

In our last example, we study a scenario with a discontinuity. We consider the same setting as in Subsection \ref{sim1/2}  except that we choose $\alpha = 1.5$ and $T = 0.008$. In this case, we suspect from Figure \ref{fig:example2} a discontinuity of the true physical solution around $t = 0.002$. 
\vspace{-13pt}
\vspace{-7pt}
The numerical approximations approach the discontinuity ``from the right'' as $\Delta\downarrow 0$, which reflects how larger time steps slow down the build-up of the probability mass close to the origin and hence delay the onset of the blow-up. By the mechanism revealed in our key Lemma \ref{lem:jump-dens}, however, for small enough $\Delta>0$, the probability of the numerical solution being close to the origin is nearly that of the physical solution just before the blow-up time of the latter, and then, by Lemma \ref{lem::special_case}, the numerical solution catches up with the discontinuity after a short time.

Figure \ref{fig:example2} also explains the choice of the M1 topology in Theorem \ref{thm::main}, rather than the J1 or the supremum norm topology (see \cite[Section 3.3]{Whi}, \cite[Appendix A.1]{CRS} for juxtapositions of these topologies), as the numerical solution seems to approach the discontinuity by many small jumps.  

We are unable to make any general claims about the rate of convergence in the blow-up regime, and therefore do not estimate it numerically. Instead, we refer the interested reader to \cite[Section 4.3]{KR} for further numerical illustrations of the convergence in this regime. 




\section*{Appendix: convergence rate of the particle approximation}

In this appendix, we specify the particle approximation to the time-stepping scheme used in Section \ref{sec::numerical}, as per \cite[Algorithm 1]{KR}, and bound the convergence rate of this approximation. Fix some $\Delta>0$, let $\{X^{\Delta,N,i}_0\}_{\stackrel{i=1,\,2,\,\ldots,\,N}{N=1,\,2,\,\ldots}}$ be i.i.d.~with the law of $X_{0-}$, and pick i.i.d.~standard Brownian motions $\{B^i\}_{i=1,\,2,\,\ldots}$. Define $X^{\Delta,N,i}_{n\Delta}$ and $\Lambda^{\Delta,N}_{n\Delta}$, $n=1,\,2,\,\ldots$ recursively~by 
\begin{align}
\label{def_part}
&\Lambda^{\Delta,N}_{n\Delta} = \frac{\alpha}{N}\sum_{j=1}^N \mathbf{1}_{\{\min_{0\le m\le n-1} X^{\Delta,N,j}_{m\Delta}<0\}}, \\
&X^{\Delta,N,i}_{n\Delta} = X^{\Delta,N,i}_0+B^i_{n\Delta}-\Lambda^{\Delta,N}_{n\Delta}.
\end{align}

\begin{prop*}
Let $\Lambda^\Delta$ be defined by \eqref{eq::init_discrete} and $\Lambda^{\Delta,N}$ by \eqref{def_part}.
For every $n_0=1,\,2,\,\ldots\,$, there exists some $C=C(\alpha,\Delta,n_0)\ge0$ such that
\begin{equation}\label{part_conv}
\E\Big[\max_{1\le n\le n_0} |\Lambda^{\Delta,N}_{n\Delta}-\Lambda^\Delta_{n\Delta}| \Big]\le\frac{C}{\sqrt{N}},\quad N=1,\,2,\,\ldots.
\end{equation}
\end{prop*}

\noindent\textbf{Proof.} With the auxiliary particle processes 
\begin{equation}
\widetilde{X}^{\Delta,N,i}_{n\Delta}:=X^{\Delta,N,i}_0+B^i_{n\Delta}-\alpha\pp\Big(\min_{0\le m\le n-1} \widetilde{X}^{\Delta,N,i}_{m\Delta}<0\Big),\quad n=1,\,2,\,\ldots,
\end{equation}
indexed by $i=1,\,2,\,\ldots,\,N$ and $N=1,\,2,\,\ldots\,$, consider the decomposition
\begin{equation*}
\begin{split}
\Lambda^{\Delta,N}_{n\Delta}-\Lambda^\Delta_{n\Delta}
=&\;\frac{\alpha}{N}\sum_{j=1}^N \mathbf{1}_{\{\min_{0\le m\le n-1} X^{\Delta,N,j}_{m\Delta}<0\}}
-\frac{\alpha}{N}\sum_{j=1}^N \mathbf{1}_{\{\min_{0\le m\le n-1} \widetilde{X}^{\Delta,N,j}_{m\Delta}<0\}} \\
&+\frac{\alpha}{N}\sum_{j=1}^N \mathbf{1}_{\{\min_{0\le m\le n-1} \widetilde{X}^{\Delta,N,j}_{m\Delta}<0\}}-\alpha\pp\Big(\min_{0\le m\le n-1} X^\Delta_{m\Delta}<0\Big).
\end{split}
\end{equation*}
Applying the triangle inequality and the upper bound $\sqrt{Np(1-p)}\le\sqrt{N}/2$ on the absolute first centered moment of a binomial distribution with parameters $N$, $p$, upon taking the absolute value and the expectation on both sides, we deduce that
\begin{equation*}
\begin{split}
\E\big[|\Lambda^{\Delta,N}_{n\Delta}-\Lambda^\Delta_{n\Delta}|\big]
\le &\;\frac{\alpha}{N}\sum_{j=1}^N \pp\Big(\min_{0\le m\le n-1} X^{\Delta,N,j}_{m\Delta}<0\le \min_{0\le m\le n-1} \widetilde{X}^{\Delta,N,j}_{m\Delta}\Big) \\
&+\frac{\alpha}{N}\sum_{j=1}^N \pp\Big(\min_{0\le m\le n-1} \widetilde{X}^{\Delta,N,j}_{m\Delta}<0\le \min_{0\le m\le n-1} X^{\Delta,N,j}_{m\Delta}\Big)+\frac{\alpha}{2\sqrt{N}}. 
\end{split}
\end{equation*} 

\smallskip

Thanks to $\min_{0\le m \le n-1} \widetilde{X}^{\Delta,N,j}_{m\Delta}\!\ge\!0\Longleftrightarrow\widetilde{X}^{\Delta,N,j}_{m\Delta}\!\ge\!0$, $0\!\le\! m\!\le\! n-1$, its analogue $\min_{0\le m \le n-1} X^{\Delta,N,j}_{m\Delta}\ge0\Longleftrightarrow X^{\Delta,N,j}_{m\Delta}\ge0$, $0\le m\le n-1$, and the union bound we have
\begin{equation*}
\begin{split}
\E\big[|\Lambda^{\Delta,N}_{n\Delta}-\Lambda^\Delta_{n\Delta}|\big]
\le &\;\frac{\alpha}{N}\sum_{j=1}^N\sum_{m=0}^{n-1} \pp(X^{\Delta,N,j}_{m\Delta}<0\le  \widetilde{X}^{\Delta,N,j}_{m\Delta}) \\
&+\frac{\alpha}{N}\sum_{j=1}^N\sum_{m=0}^{n-1} \pp(\widetilde{X}^{\Delta,N,j}_{m\Delta}<0\le X^{\Delta,N,j}_{m\Delta})+\frac{\alpha}{2\sqrt{N}}. 
\end{split}
\end{equation*} 
Noting that $X^{\Delta,N,j}_{m\Delta}=X^{\Delta,N,j}_0+B^j_{(m-1)\Delta}+(B^j_{m\Delta}-B^j_{(m-1)\Delta})
-\Lambda^{\Delta,N}_{m\Delta}$, $m=1,\,2,\,\ldots$ and $\widetilde{X}^{\Delta,N,j}_{m\Delta}=X^{\Delta,N,j}_0+B^j_{(m-1)\Delta}+(B^j_{m\Delta}-B^j_{(m-1)\Delta})
-\Lambda^\Delta_{m\Delta}$, $m=1,\,2,\,\ldots\,$, we condition on $X^{\Delta,N,j}_0+B^j_{(m-1)\Delta}
-\Lambda^{\Delta,N}_{m\Delta}$ and $X^{\Delta,N,j}_0+B^j_{(m-1)\Delta}
-\Lambda^\Delta_{m\Delta}$ in the two probabilities, and use the bound $\frac{1}{\sqrt{2\pi\Delta}}$ on the density of the thereof independent $B^j_{m\Delta}-B^j_{(m-1)\Delta}$ to arrive at
\begin{equation*}
\begin{split}
\E\big[|\Lambda^{\Delta,N}_{n\Delta}-\Lambda^\Delta_{n\Delta}|\big]
&\le\frac{\alpha}{N\sqrt{2\pi\Delta}}\sum_{j=1}^N\sum_{m=1}^{n-1} \E\big[|\Lambda^{\Delta,N}_{m\Delta}-\Lambda^\Delta_{m\Delta}|\big]+\frac{\alpha}{2\sqrt{N}} \\
&= \frac{\alpha}{\sqrt{2\pi\Delta}}\sum_{m=1}^{n-1} \E\big[|\Lambda^{\Delta,N}_{m\Delta}-\Lambda^\Delta_{m\Delta}|\big]+\frac{\alpha}{2\sqrt{N}}. 
\end{split}
\end{equation*}
It is easy to see by induction over $n=1,\,2,\,\ldots,\,n_0$ that, for some $\widehat{C}(n)=\widehat{C}(\alpha,\Delta,n)\ge0$,
\begin{equation}
\E\big[|\Lambda^{\Delta,N}_{n\Delta}-\Lambda^\Delta_{n\Delta}|\big]\le\frac{\widehat{C}(n)}{\sqrt{N}},\quad N=1,\,2,\,\ldots.
\end{equation}
We get \eqref{part_conv} with $C:=\sum_{n=1}^{n_0} \widehat{C}(n)$ via $\max_{1\le n\le n_0} |\Lambda^{\Delta,N}_{n\Delta}-\Lambda^\Delta_{n\Delta}|\le\sum_{n=1}^{n_0}  |\Lambda^{\Delta,N}_{n\Delta}-\Lambda^\Delta_{n\Delta}|$. \qed

\bigskip

\bibliography{bibliography}
\bibliographystyle{amsplain}

\bigskip

\end{document}